\documentclass[11pt]{article}
\usepackage{amsmath}
\usepackage{amsthm}
\usepackage{amssymb}
\usepackage{tikz}
\usetikzlibrary{arrows.meta}
\usepackage{subcaption} 
\usepackage{booktabs} 
\usepackage{tabularx} 

\usepackage[margin=1in]{geometry}

\usepackage[linktocpage=true]{hyperref}
\usepackage{setspace}
\usepackage{amssymb, amsmath, amsthm, graphicx,mathrsfs}
\usepackage{caption,cite}
\usepackage{mathtools,color}

\usepackage{cleveref}
\usepackage{tikz,verbatim}
\usetikzlibrary{shapes,snakes}
\usetikzlibrary{arrows}

\newtheorem{thm}{Theorem}[section]
\newtheorem{cor}[thm]{Corollary}

\newtheorem{definition}{Definition}

\newtheorem{lem}[thm]{Lemma}
\newtheorem{conj}{Conjecture}
\newtheorem{prop}[thm]{Proposition}
\newtheorem{claim}[thm]{Claim}

\newcommand{\eps}{\varepsilon}

\title{Sidorenko inequalities and symmetry-forcing\\ for graph orientations}

\author{Jiaxi Nie\thanks{School of Mathematics, Georgia Institute of Technology, 686 Cherry Street NW, Atlanta, GA 30332, USA. E-mail: {\tt jnie47@gatech.edu}} \and Fan Wei\thanks{Department of Mathematics, Duke University, 120 Science Drive, Durham, NC 27710, USA. Research supported by NSF grant  DMS-2401414.  E-mail: {\tt fan.wei@duke.edu}} }

\begin{document}

\maketitle

\begin{abstract}
We study directed graphs whose homomorphism densities never increase
or never decrease when the host edge directions are randomized. We
call these properties oriented Sidorenko (OS) and oriented
anti-Sidorenko (OAS), respectively. A directed graph is
symmetry-forcing if equality between its original and randomized
densities forces quasi-random edge directions. We characterize
connected OS directed graphs in terms of weakly norming graphs and
show that they are symmetry-forcing exactly when their underlying
graphs contain a cycle. We also prove that a connected graph is step
forcing exactly when it is step Sidorenko and contains a cycle,
giving a negative answer to a question of Dole\v{z}al et al.
Both forcing results follow from a characterization of strict
convexity: the homomorphism density functional of a graph is strictly
convex on nonnegative kernels exactly when the graph is weakly
norming and contains a cycle. On the OAS side, we construct an
infinite symmetry-forcing family and a nonbipartite example. We also
classify OS and OAS orientations of stars, paths, and cycles.
\end{abstract}

\tableofcontents

\section{Introduction}
\subsection{Quasi-random graphs, Sidorenko's conjecture and the forcing conjecture}
Since the work of Thomason~\cite{thomason1987pseudo,thomason1987random},
quasi-random graphs have become a central object in extremal
combinatorics. Let $(G_n)_{n\geq1}$ be a sequence of graphs with
$|V(G_n)|\to\infty$ and edge density $p+o(1)$. The sequence is called
\emph{quasi-random} if its edges are distributed uniformly between all
pairs of vertex sets: uniformly over $A,B\subseteq V(G_n)$,
\[
    e_{G_n}(A,B)
    =
    p|A||B|+o\bigl(|V(G_n)|^2\bigr),
\]
where $e_{G_n}(A,B)$ denotes the number of ordered pairs
$(x,y)\in A\times B$ such that $xy$ is an edge of $G_n$. This is a
global condition, since it controls the edge distribution between every
pair of vertex sets.

Chung, Graham, and Wilson~\cite{chung1989quasirandom} proved that
quasi-randomness has many equivalent descriptions. Perhaps the most
surprising one is that the global condition above is determined by the
count of a single small graph. For graphs $H$ and $G$, let $t(H,G)$ denote the \emph{homomorphism density} of $H$ in $G$, namely the probability that a uniformly random map from $V(H)$ to $V(G)$ sends every edge of $H$ to an edge of $G$. Chung, Graham, and Wilson showed that a sequence of edge density $p+o(1)$ is quasi-random if and only if $t(C_4,G_n)=p^4+o(1).$ Thus, once the edge density is fixed, quasirandomness can be characterized by the density of $C_4$.

Homomorphism densities and the special role of $C_4$ are also at the center of two major conjectures. Sidorenko's conjecture asks for a lower bound that holds in every host graph. The forcing conjecture asks what equality in that lower bound says about the host graph.

\begin{conj}[Sidorenko and forcing conjectures]\label{conj:forcing} Let $H$ be a bipartite graph. Then \begin{equation}\label{eq:Sid} t(H,G)\geq t(K_2,G)^{e(H)} \end{equation} for every graph $G$. Moreover, if $H$ contains a cycle, then for every fixed $p\in(0,1)$, any sequence $(G_n)_{n\geq1}$ with $|V(G_n)|\to\infty$ satisfying $t(K_2,G_n)=p+o(1)$ and $t(H,G_n)=p^{e(H)}+o(1)$ is quasi-random. 
\end{conj}
The first statement is Sidorenko's conjecture~\cite{sidorenko1991,simonovits1984extremal}, with an earlier supersaturation formulation due to Erd\H{o}s and Simonovits; see~\cite{erdos1983supersaturated}. The second is the forcing conjecture~\cite{skokan2004,conlon2010}.
We say a graph is Sidorenko if it satisfies~\eqref{eq:Sid}, and it is forcing if equality only holds for sequence of quasirandom graphs. The graph $C_4$ is the first and best-known example satisfying both conjectures, by the Cauchy--Schwarz inequality and the Chung--Graham--Wilson theorem. {Sidorenko's conjecture has been open for a long time, until very recently, OpenAI~\cite{OpenAI2026Math} claimed a counterexample to the Sidorenko conjecture, which, if correct, also disprove the forcing conjecture. However, the following weaker statement might still be true: every Sidorenko graph containing a cycle is forcing. We call this the \textit{Sidorenko-forcing implication}.}

\subsection{Oriented Sidorenko/anti-Sidorenko and symmetry-forcing}

Griffiths~\cite{griffiths2013quasi} initiated the study of the corresponding questions for directed graphs. Throughout, a directed graph means an orientation of a finite simple graph; in particular, loops and antiparallel edges are excluded. Here the underlying graph of the host is kept fixed, and one studies only the effect of its edge directions. One compares the density of a directed graph in a directed host with its density in a randomly oriented version of the same underlying host graph. Let $D$ and $G$ be directed graphs. The directed homomorphism density $t(D,G)$ is the probability that a uniformly random map from $V(D)$ to $V(G)$ preserves the direction of every edge. Let $\bar D, \bar G$ be the underlying graphs of $D$ and $G$ respectively. Let $G^{\mathrm{rand}}$ be obtained by keeping the underlying graph of $G$ and orienting each edge of $\bar G$ independently in either direction with probability $1/2$. Then we write  $t\left(D,\frac{G+G^T}{2}\right) := 2^{-|E(D)|}t(\bar D, \bar G) = \mathbb E\,t(D,G^{\mathrm{rand}}) + o_{|V(G)|}(1)$ to be asymptotically the expected density of $D$ after the direction of each edge in the host is randomized. 

The relevant notion of quasirandomness concerns only the edge directions. A sequence of directed graphs $(G_n)_{n\geq1}$ is \emph{quasi-random in its orientations} if \[ \max_{A,B\subseteq V(G_n)} \left| e_{G_n}(A,B)-e_{G_n}(B,A) \right| = o\left(|V(G_n)|^2\right), \] where $e_{G_n}(A,B)$ denotes the number of edges directed from $A$ to $B$. Thus no two vertex sets have a substantial bias toward one direction. The underlying undirected graphs need not be quasi-random.

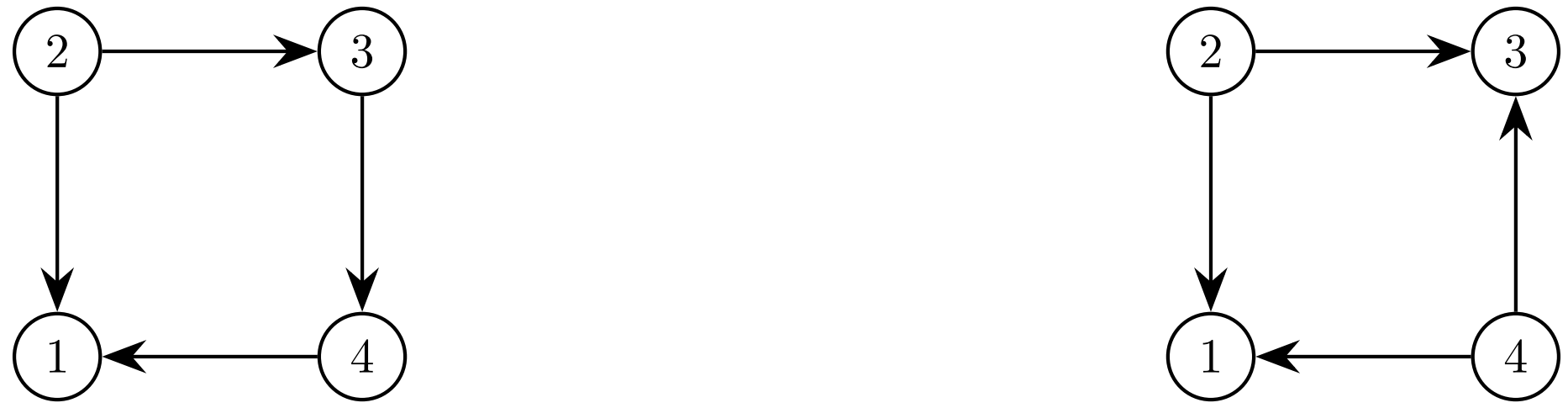
\begin{figure}[h]\centering\begin{subfigure}[b]{0.45\textwidth}\centering\begin{tikzpicture}[>={Stealth[scale=1.5]}, thick, scale=1.2]\node[circle, draw, minimum size=0.6cm] (1) at (0, 0) {1};\node[circle, draw, minimum size=0.6cm] (2) at (0, 2) {2};\node[circle, draw, minimum size=0.6cm] (3) at (2, 2) {3};\node[circle, draw, minimum size=0.6cm] (4) at (2, 0) {4};\draw[->] (2) -- (1);\draw[->] (2) -- (3);\draw[->] (3) -- (4);\draw[->] (4) -- (1);\end{tikzpicture}\caption{$O_1$}\end{subfigure}\hfill\begin{subfigure}[b]{0.45\textwidth}\centering\begin{tikzpicture}[>={Stealth[scale=1.5]}, thick, scale=1.2]\node[circle, draw, minimum size=0.6cm] (1) at (0, 0) {1};\node[circle, draw, minimum size=0.6cm] (2) at (0, 2) {2};\node[circle, draw, minimum size=0.6cm] (3) at (2, 2) {3};\node[circle, draw, minimum size=0.6cm] (4) at (2, 0) {4};\draw[->] (2) -- (1);\draw[->] (2) -- (3);\draw[->] (4) -- (3);\draw[->] (4) -- (1);\end{tikzpicture}\caption{$O_2$}\end{subfigure}\caption{Two orientations of $C_4$ studied by Griffiths. 
The orientation $O_1$ is OAS and symmetry-forcing, while $O_2$ is OS and symmetry-forcing.
}\label{fig:C4_orientations}\end{figure}

Griffiths's main purpose was to prove an analogue of Chung--Graham--Wilson theorem in which certain orientations of $C_4$ forces quasi-random orientations. He studied the two orientations $O_1$ and $O_2$ of $C_4$ shown in Figure~\ref{fig:C4_orientations}: $O_1\colon 1\leftarrow2\rightarrow3\rightarrow4\rightarrow1,$ and $O_2\colon 1\leftarrow2\rightarrow3\leftarrow4\rightarrow1. $
Griffiths proved that, for either $D=O_1$ or $D=O_2$, equality \[ t(D,G_n) = \,t\left(D,\frac{G_n+G_n^T}{2}\right)+o(1) \] holds precisely when $(G_n)_{n\geq1}$ is quasi-random in its orientations. 

Although the equality statement was the main point of his work, his proof gives more. It shows that, for every oriented host graph $G$, $t(O_1,G) \leq  t\left(O_1,\frac{G+G^T}{2}\right),$  while $t(O_2,G) \geq t\left(O_2,\frac{G+G^T}{2}\right)$.  Thus the two orientations satisfy opposite inequalities for every host, not only a statement about the equality case.

We introduce names for these two properties. 
\begin{definition}[Oriented Sidorenko (OS) and oriented anti-Sidorenko (OAS)] 
A directed graph $D$ has the \emph{oriented Sidorenko property}, or is \emph{OS}, if \[ t(D,G) \geq  t\left(D,\frac{G+G^T}{2}\right) \] for every directed graph $G$. A directed graph $D$ has the \emph{oriented anti-Sidorenko property}, or is \emph{OAS}, if \[ t(D,G) \leq t\left(D,\frac{G+G^T}{2}\right) \] for every directed graph $G$. \end{definition}

Equivalently, these properties may be defined using digraphons. We give
the continuous formulation in Section~\ref{sec:preliminaries}, where it
will be used in the proofs.

We also give a name to the equality property studied by Griffiths. \begin{definition}[Symmetry-forcing]\label{def:symmetry-forcing} A directed graph $D$ is \emph{symmetry-forcing} if, for every sequence of directed graphs $(G_n)_{n\geq1}$ with $|V(G_n)|\to\infty$, \begin{equation}\label{eq:symmetry-forcing}
    t(D,{G}_n)
    =
    t\left(D,\frac{{G}_n+{G}_n^T}{2}\right)+o(1)    
\end{equation} 
holds if and only if $(G_n)_{n\geq1}$ is quasi-random in its orientations. \end{definition}

 Equivalently, in the language of digraphons, $D$ is symmetry-forcing if, for every digraphon $W$, the equality $t(D,W)=t\left(D,\frac{W+W^T}{2}\right)$
 holds if and only if $W=W^T$ almost everywhere. (see Section~\ref{sec:preliminaries} for definitions of digraphon and its density)

In this terminology, Griffiths proved the following. 
\begin{thm}[Griffiths~\cite{griffiths2013quasi}] \label{prop:C4} The directed graph $O_1$ is OAS and symmetry-forcing, while $O_2$ is OS and symmetry-forcing.
\end{thm} 
The OS and OAS properties concern a universal inequality. Symmetry forcing concerns the equality case. This is the same distinction that appears between Sidorenko's conjecture and the forcing conjecture. There is also a new feature in the orientation setting: the universal inequality can go in either direction.

Griffiths's result therefore does more than give two tests for quasi-random orientations. It gives the first examples of two classes of directed graphs defined by opposite universal inequalities. This leads to the problem studied in this paper: which directed graphs are OS, which are OAS, and what structure is forced by these properties? We give a complete characterization for connected OS directed graphs in terms of weakly norming graphs, and further show that a connected OS directed graph is symmetry forcing if and only if its underlying graph contains a cycle. On the OAS side, we construct an infinite family of symmetry-forcing examples and give an example whose underlying graph is not bipartite. Moreover, we characterize all OS and OAS orientations of stars, paths, and cycles. We also prove that every symmetry-forcing directed graph must be either OS or OAS. For brevity, we have stated the definitions above using finite directed graphs. The equivalent digraphon formulations are given in Section~\ref{sec:preliminaries} and are used in the proofs.

 We briefly explain why focusing on   ``quasi-random orientations", rather than a fully uniform directed graph process like ${G}(n,p)$, is both combinatorially necessary and richer.
 First, it prevents the rigid skeleton structure (the set of underlying edges in $ G$) of a host graph from obscuring a perfectly random flow. For example, as observed in \cite{griffiths2013quasi}, randomly oriented bipartite graph has maximum directional entropy, but would trivially fail a standard ${G}(n,p)$ test simply because it lacks triangles. More importantly, decoupling the orientation is precisely what unlocks the Oriented anti-Sidorenko (OAS) property. In the classical undirected setting, absolute anti-Sidorenko bounds are structurally impossible for all connected graphs except a single edge, as one can always artificially inflate a subgraph's density far above the uniform baseline simply by concentrating the host graph's edges into a dense clique. By fixing the underlying skeleton and randomizing only the orientations, we remove the structural freedom needed to artificially inflate subgraph densities. This isolates the effect of the edge directions, providing a natural setting where both OS and OAS properties can exist.

\subsection{Step Sidorenko property}
{
The OS and OAS inequalities concern how homomorphism densities
change when edge directions in the host graph are averaged while the underlying host is fixed. A different averaging question has been studied in the
undirected graph setting, where averaging takes place over vertex classes
and may change the underlying host itself. This is the setting of
the \emph{step Sidorenko property}, introduced by Kr\'al', Martins,
Pach, and Wrochna~\cite{KralMartinsPachWrochna2019}. Given a graphon $U$
and a finite measurable partition $\mathcal P$ of its vertex space,
let $U_{\mathcal P}$ be obtained by averaging $U$ on each rectangle
of $\mathcal P\times\mathcal P$. A graph $H$ is step Sidorenko if
\[
    t(H,U_{\mathcal P})\leq t(H,U)
\]
for every (symmetric) graphon $U$ and vertex partition  $\mathcal P$. Further, it is \emph{step forcing} if the minimum occur only when the host are already uniform relative to the partition.  Dole\v{z}al, Greb{\'\i}k, Hladk\'y,
Rocha, and Rozho\v{n}~\cite{DolezalGrebikHladkyRochaRozhon2022} showed, among other things, that step Sidorenko is equivalent to \emph{cut distance compatibility}, i.e., monotonicity of $t(H,\cdot)$ under the graphon structuredness order. We defer the formal definitions to Section~\ref{sec:preliminaries}.

Our averaging operation is considerably more specific. In digraphon language, the symmetrization
\[
    W\longmapsto \frac{W+W^T}{2} 
\]
is conditional expectation onto the $\sigma$-algebra invariant under the coordinate swap $(x,y)\mapsto(y,x)$: it retains each unordered pair and averages only its two possible directions. This is not, in general, averaging over rectangles induced by a partition of the vertex space. Nevertheless, there is a concrete specialization of step averaging that realizes our operation. If $W$ is a finite step digraphon, form the equal-measure disjoint union $B_W=W\oplus W^T$ and let $\mathcal Q_W$ be the specific vertex partition merging corresponding step parts of the two copies. Then $(B_W)_{\mathcal Q_W}$ is, up to measure-preserving relabelling, $(W+W^T)/4$. For a connected directed graph satisfying $D\cong D^T$, a tensor-power argument shows that $D$ is OS if and only if
\[
    t\bigl(D,(B_W)_{\mathcal Q_W}\bigr)\leq t(D,B_W)
\]
for every such $W$; see Claim~\ref{prop:OS-paired-stepping}.
Thus OS can be viewed as a specialization of the directed analogue of step Sidorenko to these particular doubled hosts and paired partitions; while OAS concerns the reverse inequality under the orientation symmetrization itself.}

\subsection{Our Results}

\subsubsection{OS directed graphs}

Our first main theorem provides a complete characterization of connected
OS directed graphs, i.e., a complete characterization of the
connected directed graphs for which symmetrization cannot increase the
homomorphism density. 

For a directed graph $D$, let $D^T$ be obtained by reversing every edge, and let $\overline D$ be its underlying undirected graph.

 \begin{thm}\label{thm:charcOS} Let $D$ be a connected directed graph with at least one edge. Then $D$ is OS if and only if all of the following conditions hold: \begin{itemize} \item[(a)] $D$ is a {\it one-direction orientation}, i.e., its vertex set has a bipartition $A\cup B$ such that every edge is directed from $A$ to $B$; 
 \item[(b)] $D$ is isomorphic to $D^T$; \item[(c)] the underlying graph $\overline D$ is weakly norming. \end{itemize} 
 \end{thm}

This is an exact characterization. Every connected OS directed graph must have all three properties, and the three properties are also sufficient. In particular, every connected OS directed graph is bipartite, and every vertex is either a source or a sink. Reversing all edges must give the same directed graph up to isomorphism. We defere the definitions of (weakly) norming to Section~\ref{sec:preliminaries}.

Theorem~\ref{thm:charcOS} gives many examples of OS directed graphs. These include the one-direction orientations of even cycles, balanced complete bipartite graphs $K_{t,t}$, $K_{t,t}$ minus a perfect matching, and hypercubes\cite{hatami2010,lee2021convex}. 

{
Our second main result characterizes equality in the OS inequality. Recall that the graph Sidorenko-forcing implication says every Sidorenko graph containing a cycle is forcing. We prove an OS version of the Sidorenko-forcing implication.

\begin{thm}\label{thm:OS_forcing}
A connected OS directed graph is symmetry-forcing if and only if its underlying graph contains a cycle.   
\end{thm}

Together with Theorem~\ref{thm:charcOS}, this gives an exact characterization of connected OS symmetry-forcing directed graphs.
}

\subsubsection{Consequences in the undirected setting}

Kr\'al', Martins, Pach, and Wrochna
~\cite{KralMartinsPachWrochna2019} observed, using
Lov\'asz~\cite[Proposition~14.13]{lovasz2012large}, that every weakly
norming graph has the step Sidorenko property, and asked whether the
converse holds for connected graphs. Dole\v{z}al et
al.~\cite{DolezalGrebikHladkyRochaRozhon2022} proved the converse, i.e. every connected graph with the step Sidorenko property is weakly norming.
Consequently, combining their result with Theorem~\ref{thm:charcOS}, for a connected one-direction orientation
$D\cong D^T$ with at least one edge we have
\[
    D\text{ is OS}
    \quad\Longleftrightarrow\quad
    \overline D\text{ is step Sidorenko}.
\]

{
For the equality case, Dole\v{z}al et al.~\cite{DolezalGrebikHladkyRochaRozhon2022} proved that if a connected undirected graph $H$ is norming, then the graphon parameter $t(H,\cdot)$ is \emph{cut distance identifying}, and hence $H$ is step forcing. They asked whether the converse is true~\cite[Problem 3.29]{DolezalGrebikHladkyRochaRozhon2022}, i.e. does undirected graph $H$ being step forcing (or $t(H,\cdot)$ being {cut distance identifying}) imply that $H$ is norming? We defer the formal definitions to Section~\ref{sec:preliminaries}. 

The following theorem answers~\cite[Problem 3.29]{DolezalGrebikHladkyRochaRozhon2022} in the negative.

\begin{thm}\label{thm:cyclic-step-forcing}
Let $H$ be a weakly norming graph containing a cycle. Then $t(H,\cdot)$ is a cut distance identifying graphon parameter, and hence $H$ is step forcing.
\end{thm}
Indeed, some weakly norming graphs containing a cycle, such as $K_{2,3}$, are not norming (see Table~\ref{tab:weakly-norming-not-norming}), whereas the theorem shows that all such graphs are step forcing. Thus, step forcing does not imply norming. 

\begin{table}[h]
    \centering
    \small
    \renewcommand{\arraystretch}{1.3}
    \begin{tabularx}{\textwidth}{@{}XX@{}}
        \toprule
        Family & Reason or reference \\
        \midrule
        Complete bipartite graphs \(K_{m,n}\), where \(m,n\ge 2\)
        and at least one of \(m,n\) is odd
        & Odd-degree obstruction: every norming graph has
        even vertex degrees.~\cite{hatami2010} \\

        Hypercubes \(Q_d\), where \(d\ge 3\)
        & Lee--Sidorenko~\cite{LeeSidorenko2022}, Theorem~1.5. \\

        Crown graphs \(K_{n,n}-\text{perfect matching}\), where \(n\ge 4\)
        & Lee--Sch\"ulke~\cite{lee2021convex}, Theorem~1.6. \\

        \(1\)-subdivisions\footnotemark{}  of complete graphs $K_n$,
        where \(n\ge 4\)
        & Lee--Sidorenko~\cite{LeeSidorenko2022},  Theorem~6.6. \\
        \bottomrule
    \end{tabularx}
    \caption{Representative families of weakly norming graphs that are not norming.}
    \label{tab:weakly-norming-not-norming}
\end{table}
\footnotetext{The \(1\)-subdivision of a graph is obtained
by replacing every edge by a path of length two.}

Theorem~\ref{thm:cyclic-step-forcing} together with the result of Doležal et al. that a connected graph is step Sidorenko if and only if it is weakly norming gives the following corollary, which can be viewed as the step version of the Sidorenko-forcing implication.
\begin{cor}\label{cor:char-step-forcing}
A connected graph is step forcing if and only if it is step Sidorenko and contains a cycle.
\end{cor}

Lee and Schülke~\cite{lee2021convex} characterized weakly norming graphs by the convexity of their homomorphism density functionals on graphons, and norming graphs by strict convexity on signed graphons. Their proofs also adapt to the corresponding bipartite homomorphism density functionals on possibly asymmetric kernels.

Strict convexity on nonnegative kernels, however, appears to be less well understood. We prove the following characterization. For a bipartite graph $H$ with a fixed bipartition $X\cup Y$ and any kernel $W$, $t_{bi}(H,W)$ is the total weight of homomorphisms from $H$ to $W$ where vertices in $X$ are mapped only to the first coordinate space and vertices in $Y$ are mapped only to the second coordinate space. See Section~\ref{sec:preliminaries} for formal definitions.

\begin{thm}\label{thm:strict-convexity}
Let $H=(X\cup Y,E)$ be a bipartite graph with a fixed bipartition and at least one edge. Then $H$ is weakly norming and contains a cycle if and only if $t_{bi}(H,\cdot)$ is strictly convex on nonnegative kernels.
\end{thm}
}

\subsubsection{OAS Directed Graphs}

A complete characterization of OAS directed graphs remains open. We prove several results that show that this side has a different structure.

 We have shown there are infinitely many OS symmetry-forcing graphs. In this subsection, we also construct an infinite family $\{D_t\}_{t\ge 2}$ (see Definition \ref{def:Dt}) of OAS symmetry-forcing directed graphs, and find an example of an OAS directed graph whose underlying graph is not bipartite. 
 
 For every $t\geq2$, let $D_t$ be the orientation of $K_{t,2^t}$ in which the vertices on the side of size $2^t$ realize all possible orientations of their incident edges to the side of size $t$. The precise definition is given in Section~\ref{sec:OAS}.

\begin{thm}\label{thm:Dt OAS symmetry}
For every $t\geq 2$, the directed graph $D_t$ is OAS and
symmetry-forcing. 
\end{thm}

The non-bipartite OAS example has an underlying graph containing a triangle. This gives a clear difference between OAS and OS. By Theorem~\ref{thm:charcOS}, every connected OS directed graph has a bipartite underlying graph, but the underlying graph of an OAS directed graph need not be bipartite.

\begin{prop}\label{prop:nonbipartite-oas}
There exists an OAS directed graph whose underlying graph contains a
triangle.
\end{prop}

{
This also gives an OAS directed graph whose underlying graph contains
a cycle but which is not symmetry-forcing. Indeed, let $W$ be the
digraphon corresponding to a single directed edge. If the underlying
graph of $D$ is not bipartite, then $t(D,W)=t\left(D,\frac{W+W^T}{2}\right)=0$ while $W\ne W^T$. Thus the OAS analogue of
Theorem~\ref{thm:OS_forcing} fails. It remains open to characterize
connected OAS symmetry-forcing directed graphs, and this problem
appears substantially more complicated than the OS case.
}

\subsubsection{Complete characterization of OS/OAS stars, paths, and cycles}
For stars, paths, and cycles, we obtain a complete characterization.

\begin{thm}\label{thm:OAS stars}
Let $D$ be an oriented star.
\begin{enumerate}
    \item If $D$ has one edge, then $D$ is both OS and OAS.
    \item If $D$ has at least two edges, then $D$ is not OS. Moreover,
    $D$ is OAS if and only if its center has equal in-degree and
    out-degree.
\end{enumerate}
\end{thm}

\begin{thm}\label{thm:OAS paths}
Let $P$ be an oriented path. 
\begin{enumerate}
    \item If $P$ has one edge, then $P$ is both OS and OAS.
    \item If $P$ has at least two edges, then $P$ is not OS. Moreover,
    $P$ is OAS if and only if it is the directed two-edge path $1\rightarrow 2\rightarrow 3$.
\end{enumerate}
\end{thm}

\begin{thm}\label{thm:OAS cycles}
An oriented cycle is OS if and only if it is an alternating orientation of an even cycle, and it is OAS if and only if it is isomorphic to $O_1$.
\end{thm}

Our characterization of connected OS directed graphs (Theorem \ref{thm:charcOS}) already determines all OS orientations of cycles: they are exactly the even cycles whose vertices alternate between sources and sinks. It therefore remains only to classify the OAS orientations of cycles. We prove that, up to isomorphism, $O_1$ is the only one. Together, these results give a complete classification of all oriented cycles that are OS or OAS.

\subsubsection{The relation with symmetry forcing} We return to the equality question. Although OS, OAS, and symmetry forcing are different properties, every symmetry-forcing graph must satisfy one of the two universal inequalities. \begin{prop}\label{prop:forcingandOS} If a directed graph $D$ is symmetry-forcing, then $D$ is either OS or OAS. \end{prop} Thus a symmetry-forcing graph always has one fixed inequality direction. We have characterized connected OS symmetry-forcing directed graphs. A corresponding characterization of connected OAS symmetry-forcing directed graphs remains open.

\subsection{Relation with other directed Sidorenko properties}
Another analogue of the Sidorenko property was recently introduced by Fox, Himwich, Mani, and Zhou \cite{fox2025directed}. Given a fixed directed graph \(F\), they say that \(F\) has the \emph{directed Sidorenko (DS) property} if, for every digraphon \(W:[0,1]^2\to[0,1]\), the homomorphism density of \(F\) satisfies $t(F,W)\ge t(F,pI)=p^{\,e(F)}$, where \(p=t(\vec{K}_2,W)\) and $I$ is the all-1 digraphon. Within this framework, they also isolated the important special case of tournaments~\cite{fox2024variations}: a directed graph \(F\) is said to have the \emph{tournament Sidorenko (TS) property} if the same inequality holds for every tournamenton \(W\) satisfying \(W(x,y)+W(y,x)=1\). It is interesting that there exist directed graphs $F$ (e.g., directed paths) with the so-called \emph{tournament anti-Sidorenko (TAS) property}, i.e., for every tournamenton $W$,
$t(F,W)\;\le\; t(F,pI)=p^{\,e(F)}$ where  \(p=t(\vec{K}_2,W)\). These properties have attracted considerable recent attention~\cite{fox2024variations,GRZESIK2023117,he2025new,chen2026tournament,fang2026maximizing}. One may similarly define the DAS property, but it is not as interesting since it is not hard to see that any digraph that is not a matching is not DAS.

Clearly, OS/OAS implies TS/TAS, but the relation between OS/OAS and DS/DAS is less obvious. Fox, Himwich, Mani, and Zhou \cite{fox2025directed} showed that a directed graph is directed Sidorenko if it is homomorphic to a directed edge and its underlying graph has the asymmetric-Sidorenko property (which is stronger than and conjectured to be equivalent to the Sidorenko property).  Thus Theorem~\ref{thm:charcOS} implies that, for connected directed graphs, OS is strictly stronger than DS, and hence TS.  Fox, Himwich, Mani, and Zhou characterize oriented stars in terms of the TAS property: they show that an oriented star is TAS if and only if its center's in-degree and out-degree differ by at most one. Thus, by Theorem~\ref{thm:OAS stars}, we know that OAS is strictly stronger than TAS.

Overall, for connected directed graphs, we have OS strictly stronger than DS strictly stronger than TS, and DAS strictly stronger than OAS strictly stronger than TAS.

\subsection{Structure of the Paper}
The remainder of the paper is organized as follows. In Section~\ref{sec:preliminaries}, we introduce the necessary notation and definitions. In Section~\ref{sec:forcing implies OS/OAS}, we show that every symmetry-forcing directed graph is either OS or OAS (Proposition~\ref{prop:forcingandOS}). In Section~\ref{section:OS}, we first characterize connected OS directed graphs (Theorem~\ref{thm:charcOS}), then characterize undirected bipartite graphs $H$ such that $t_{bi}(H,\cdot)$ is strictly convex on nonnegative kernels (Theorem~\ref{thm:strict-convexity}), and then use this characterization to prove the oriented version (Theorem~\ref{thm:OS_forcing}) and step version (Theorem~\ref{thm:cyclic-step-forcing}) of the Sidorenko-forcing implication. In Section~\ref{sec:OAS} we construct a family of OAS symmetry-forcing directed graphs (Theorem~\ref{thm:Dt OAS symmetry}), and an OAS directed graph with non-bipartite underlying graph (Proposition~\ref{prop:nonbipartite-oas}). Finally, in Section~\ref{sec:stars, paths, and cycles}, we characterize stars (Theorem~\ref{thm:OAS stars}), paths (Theorem~\ref{thm:OAS paths}), and cycles (Theorem~\ref{thm:OAS cycles}).

\section{Definitions, Notation, and Preliminaries}
\label{sec:preliminaries}
In this section we review some of the basic definitions and notation used in this paper. 

\medskip

\noindent\textbf{Directed graphs.}
Throughout this paper, a directed graph is an orientation of a finite simple graph; loops and antiparallel edges are excluded. For a directed graph $D$, we denote its vertex set by $V(D)$ and its directed edge (or arc) set by $E(D)$. The number of vertices and edges are denoted by $v(D) = |V(D)|$ and $e(D) = |E(D)|$, respectively. A directed edge $e$ from a vertex $u$ to a vertex $v$ is denoted as the ordered pair $(u,v)$, and we sometimes use the notation $e_+=v$ and $e_-=u$.

\medskip

\noindent\textbf{Homomorphisms and Discrete Homomorphism Densities.}
Let $H$ and $G$ be directed graphs. A \emph{homomorphism} from $H$ to $G$ is a mapping $f \colon V(H) \to V(G)$ that preserves adjacency and edge orientation; that is, if $(u,v) \in E(H)$, then $(f(u), f(v)) \in E(G)$. 

The \emph{homomorphism density} of $H$ in $G$, denoted by $t(H, G)$, represents the probability that a uniformly random mapping from $V(H)$ to $V(G)$ is a valid homomorphism. Formally, it is defined as $t(H, G) = \frac{|\mathrm{Hom}(H, G)|}{|V(G)|^{v(H)}}$, where $\mathrm{Hom}(H, G)$ denotes the set of all homomorphisms from $H$ to $G$. 

\medskip

\noindent\textbf{Graphons, Digraphons and Continuous Homomorphism Densities}
To study the asymptotic properties of dense graph sequences, it is natural to extend the discrete homomorphism density to the continuous setting using the theory of graph limits. 

A (directed) \emph{kernel} is a bounded measurable function $U\colon[0,1]^2\to\mathbb R$, not necessarily symmetric. We identify kernels that agree almost everywhere. A \emph{digraphon} is a kernel with codomain restricted to $[0,1]$. Its transpose is $W^T(x,y)=W(y,x)$, and its symmetrization is $S=\frac{W+W^T}{2}$. A \emph{graphon} is just a symmetric digraphon, i.e. it satisfies $W(x,y)=W(y,x)$ for almost all $(x,y)\in [0,1]^2$.

For a fixed directed graph $H$ and a kernel $W$, the continuous \emph{homomorphism density} of $H$ in $W$ is defined by the integral
\begin{equation}\label{eq:continuous_density}
    t(H, W) = \int_{[0,1]^{V(H)}} \prod_{(u,w) \in E(H)} W(x_u, x_w) \prod_{v \in V(H)} \mathrm{d}x_v.
\end{equation}
Notice here that $(u,w)\in E(H)$ denotes the directed edge from $u$ to $w$ in $H$. 
 When $W$ is a symmetric real-valued kernel and $H$ is an undirected graph, $t(H,W)$ is defined by the same integral, independently of how the endpoints of each edge are ordered. For brevity, we sometimes omit ``$\prod_{v \in V(H)} \mathrm{d}x_v$'' when it is clear from the context, and we write $x_e:=(x_u,x_w)$ for $e=(u,w)$, i.e. $t(H, W) = \int\prod_{e\in E(H)} W(x_e)$.

Note that the discrete and continuous notions perfectly align: for any directed graph $G$, $t(H, G) = t(H, W_G)$, where $W_G$ is the step-function digraphon associated with $G$. 

\medskip

\noindent\textbf{Bipartite Homomorphism Density, Norming, and Weakly Norming} 
Let $H = (A \cup B, E)$ be a finite simple bipartite graph with at least one edge and a fixed bipartition. For any kernel $W$, the \emph{bipartite homomorphism density} of $H$ into $W$ is defined as
\begin{equation}\label{eq:bi-hom-density}
    t_{bi}(H, W) = \int_{[0,1]^{A \cup B}} \prod_{uv \in E, \, u \in A, \, v \in B} W(x_u, y_v) \prod_{u \in A} \mathrm{d}x_u \prod_{v \in B} \mathrm{d}y_v.
\end{equation}
Define $\|W\|_H:=|t_{bi}(H,W)|^{1/e(H)}$ and
$\|W\|_{r(H)}:=t_{bi}(H,|W|)^{1/e(H)}$. We say that $H$ is
\emph{norming} if $\|\cdot\|_H$ is a norm, and \emph{weakly norming}
if $\|\cdot\|_{r(H)}$ is a norm on the space of all kernels. Thus the
bipartition is fixed, and symmetry of a kernel is not assumed. We use
``norming'' for a genuine norm, as in~\cite{LeeSidorenko2022}; a graph
whose functional is only a seminorm is called \emph{seminorming}.

When we refer explicitly to the symmetric convention, these functionals
are instead restricted to symmetric kernels and $t(H,W)$ is used.
Lee and Sidorenko~\cite[Theorem~1.3]{LeeSidorenko2022} proved that
the symmetric and asymmetric conventions agree for connected bipartite
graphs. For disconnected graphs the conventions can differ, and the
choices of bipartition on different components must be retained. Unless
the symmetric convention is explicitly specified, weak norming of a
bipartite graph with a fixed bipartition always refers to the full kernel
space in this paper.

We may replace $[0,1]$ by any standard atomless probability space.
For bipartite densities, measure-preserving changes of variables may
be made independently in the two coordinates. In particular, we will
use $\Omega=[0,1]\times\{-1,1\}$, equipped with Lebesgue measure and
the uniform probability measure on the sign coordinate. We shall also
use continuity in $L^1$ on uniformly bounded kernels: if $|U|,|V|\le M$
and $m=e(H)$, telescoping the edge product gives
\begin{equation}\label{eq:kernel-L1-continuity}
    |t_{bi}(H,U)-t_{bi}(H,V)|
    \le mM^{m-1}\|U-V\|_{L^1}.
\end{equation}
The same estimate holds for directed densities and for ordinary
homomorphism densities of symmetric kernels.

We shall use the standard fact that the absolute value is redundant for even cycles; see Hatami~\cite[Section~2.3]{hatami2010}. Namely, for every real kernel $K$ and every $q\ge2$,
$$
    \|K\|_{C_{2q}}
    :=t_{bi}(C_{2q},K)^{1/(2q)}
$$
is a norm (the Schatten--von Neumann $2q$-norm). In particular, $t_{bi}(C_{2q},K)\ge0$, with equality if and only if $K=0$ almost everywhere.

\medskip

\noindent\textbf{Step Sidorenko, step forcing, cut distance compatible and cut distance identifying.}
For a graphon $W$ and a finite measurable partition $\mathcal P$ of $[0,1]$ into parts of positive measure, let $W^{\mathcal P}=W_{\mathcal P}$ be the graphon whose value on $P_i\times P_j$ is
$$
    W^{\mathcal P}(x,y)
    =\frac{1}{|P_i||P_j|}\int_{P_i\times P_j}W(s,t)\,ds\,dt
    \qquad (x\in P_i,\ y\in P_j).
$$
Following~\cite{KralMartinsPachWrochna2019}, a graph $H$ is \emph{step forcing} if
$$
    t(H,W^{\mathcal P})<t(H,W)
    \quad\text{whenever }W\text{ and }W^{\mathcal P}\text{ differ on a set of positive measure}.
$$
The non-strict version of this inequality is the \emph{step Sidorenko property}. In particular, weakly norming graphs have the step Sidorenko property~\cite{KralMartinsPachWrochna2019,DolezalGrebikHladkyRochaRozhon2022}.

We also recall the stronger notion from~\cite{DolezalGrebikHladkyRochaRozhon2022}. A sequence of graphons $(W_n)_{n\ge 1}$ weak$^*$ converges to a graphon $W$ if
$$
\int_{A\times B}W_n(x,y)dxdy\xrightarrow{n\rightarrow\infty} \int_{A\times B}W(x,y)dxdy
$$
for every measurable $A,B\subseteq[0,1]$. For a measure-preserving bijection $\phi$ of $[0,1]$, put $W^\phi(x,y)=W(\phi(x),\phi(y))$. The \emph{envelope} $\langle W\rangle$ is the set of weak$^*$ limits of sequences of versions $W^{\phi_n}$. Write $U\preceq W$ when $\langle U\rangle\subseteq\langle W\rangle$, and $U\prec W$ when the inclusion is strict. A graphon parameter $\theta$ is \emph{cut distance compatible} if $U\preceq W$ implies $\theta(U)\le \theta(W)$, and is \emph{cut distance identifying} if $U\prec W$ implies $\theta(U)<\theta(W)$. We use the cut norm $\|F\|_\square:=\sup_{S,T\subseteq[0,1]}\bigl|\int_{S\times T}F\bigr|$ and cut distance $\delta_\square(U,W):=\inf_\phi\|U-W^\phi\|_\square$. In particular, $\delta_\square(U,W)=0$ if and only if $\langle U\rangle=\langle W\rangle$~\cite{DolezalGrebikHladkyRochaRozhon2022}.

\medskip

\noindent\textbf{Oriented Sidorenko (OS), Oriented Anti-Sidorenko (OAS) and Symmetry Forcing in Digraphon Language} A directed graph $D$ is OS, respectively OAS, if for every digraphon $W$, $t(D,W)\ge t(D,\frac{W+W^T}{2})$, respectively $t(D,W)\le t(D,\frac{W+W^T}{2})$. Further, $D$ is symmetry-forcing if $t(D,W)=t(D,\frac{W+W^T}{2})$ implies $W=W^T$.

{
We also record a formulation of OS that will be useful for comparison with the step Sidorenko property. For a connected directed graph $D$ satisfying $D\cong D^T$, the following claim expresses OS as the restriction of the directed step Sidorenko inequality to a particular family of doubled hosts and paired partitions. Thus OS corresponds to a restricted subfamily of step Sidorenko tests in the directed setting.

\begin{claim}
\label{prop:OS-paired-stepping}
Let $D$ be a connected directed graph with $D\cong D^T$. For a finite step digraphon $W$, let $W\oplus W^T$ be the equal-measure disjoint union, with zero off-diagonal blocks, and let $\mathcal Q_W$ merge corresponding step parts of the two copies. Then $D$ is OS if and only if
\[
    t\bigl(D,(W\oplus W^T)_{\mathcal Q_W}\bigr)\leq t(D,W\oplus W^T)
\]
for every finite step digraphon $W$.
\end{claim}
For comparison, recall that an undirected graph $H$ has the step Sidorenko property if, for {\it every}  (symmetric) graphon $U$ and {\it every} finite measurable vertex partition $\mathcal P$, 
    $t(H,U_{\mathcal P})\leq t(H,U)$. Thus Claim~\ref{prop:OS-paired-stepping} requires the analogous inequality  for the particular doubled digraphons of forms $W\oplus W^T$ and the particular paired partitions $\mathcal Q_W$.

\begin{proof}
Write $n=v(D)$, $m=e(D)$, $S(W)=(W+W^T)/2$, and
$C=2^{m+1-n}\geq1$. Stepping the doubled host gives $S(W)/2$,
up to measure-preserving relabelling. Thus, $t\bigl(D,(W\oplus W^T)_{\mathcal Q_W}\bigr)=2^{-m}t\left(D,S(W)\right)$. Connectedness and
$D\cong D^T$ give $t(D,W\oplus W^T)=2^{1-n}t(D,W)$.
Thus the stated inequality is equivalent to
$t(D,S(W))\leq C\,t(D,W)$, which follows immediately from OS.

Conversely, assume this inequality for every finite step digraphon.
For $F\subseteq E(D)$, let $D_F$ be obtained by reversing the edges
in $F$. Expanding gives
$t(D,S(W))=2^{-m}\sum_{F\subseteq E(D)}t(D_F,W)$.
The tensor power $W^{\otimes k}$ is again a finite step digraphon
on a product probability space, which we identify with $[0,1]$.
Jensen's inequality, tensor-product multiplicativity, and the
assumed inequality imply
\[
 t(D,S(W))^k
 \leq 2^{-m}\sum_{F\subseteq E(D)}t(D_F,W)^k
 =t(D,S(W^{\otimes k}))
 \leq C\,t(D,W)^k.
\]
Taking $k$th roots and letting $k\to\infty$ proves the OS inequality
for finite step digraphons. Approximation in $L^1$ extends it to
all digraphons.
\end{proof}

}

\noindent\textbf{Decorated densities and weak H\"older inequality.} Let $H=(X\cup Y,E)$ be a bipartite graph with a fixed bipartition. A \emph{decoration} of a graph $H$ is a family $\mathbf U=(U_e)_{e\in E(H)}$ where each $U_e$ is a kernel. The \emph{decorated homomorphism density} of $H$ in $\mathbf{U}$ is defined as 
\begin{equation}
   t_{bi}(H;\mathbf U):=\int
\prod_{e=xy\in E(H)}U_e(x_x,y_y)
\prod_{x\in X}dx_x\prod_{y\in Y}dy_y. \label{eq:decorated} 
\end{equation}

Intuitively, we may interpret $\mathbf{U}$ as a collection of bipartite graphs, each representing a color, sharing the same vertex bipartition, and $t_{bi}(H;\mathbf U)$ measures the density of \emph{rainbow} homomorphisms of $H$ in $\mathbf U$ that respect the vertex bipartition and the decoration.

We shall use the weak H\"older inequality for weakly norming graphs~\cite{hatami2010}.

\begin{lem}[Weak H\"older inequality ~\cite{hatami2010}]\label{lem:weak-holder}
Let $H=(X\cup Y,E)$ be weakly norming. Let $\mathbf U=(U_e)_{e\in E(H)}$ be a decoration of $H$, then
$$
\left|
t_{bi}(H;\mathbf{U})
\right|
\le
\prod_{e\in E(H)}t_{bi}(H,|U_e|)^{1/e(H)}.
$$
\end{lem}

\medskip
\noindent\textbf{Convexity and graph norms on the full kernel space.}
Lee and Sch\"ulke~\cite[Theorem~1.5]{lee2021convex} state their
convexity characterizations for symmetric kernels. We record the
corresponding statements for our bipartite convention, together with
the short adaptation of their proofs.

\begin{lem}\label{lem:kernel-convexity}
Let $H=(X\cup Y,E)$ be a bipartite graph with a fixed bipartition and
$m=e(H)>0$.
\begin{enumerate}
    \item $H$ is weakly norming if and only if $t_{bi}(H,\cdot)$ is
    convex on nonnegative kernels.
    \item $H$ is norming if and only if $t_{bi}(H,\cdot)$ is strictly
    convex on real-valued kernels.
\end{enumerate}
\end{lem}

\begin{proof}
Write $F(W)=t_{bi}(H,W)$. If $H$ is weakly norming, then
$F(W)=\|W\|_{r(H)}^m$ on nonnegative kernels, so $F$ is convex there.
Conversely, suppose that $F$ is convex on nonnegative kernels, and put
$N(W)=F(|W|)^{1/m}$. This functional is clearly homogeneous and
positive definite. For the latter assertion, if $K$ is nonzero,
we may scale it so that $|K|\le1$. Writing $a=|X|$ and $b=|Y|$, two
applications of Jensen's inequality give
\[
    N(K)^m=F(|K|)\ge t_{bi}(K_{a,b},|K|)
    \ge\left(\int |K|\right)^{ab}>0.
\]
For nonzero kernels $U,V$, set $p=N(U)$ and $q=N(V)$. Monotonicity
on nonnegative kernels and convexity imply
\[
\begin{aligned}
    F\left(\frac{|U+V|}{p+q}\right)
    &\le F\left(\frac{p}{p+q}\frac{|U|}{p}
                 +\frac{q}{p+q}\frac{|V|}{q}\right)\\
    &\le\frac{p}{p+q}F(|U|/p)
          +\frac{q}{p+q}F(|V|/q)=1.
\end{aligned}
\]
Thus $N(U+V)\le N(U)+N(V)$, also trivially when one kernel is zero.
This proves the first statement.

For the second statement, suppose first that $H$ is norming. Every
vertex then has even degree: otherwise the nonzero kernel $f(x)f(y)$,
where $f$ takes the values $1,-1$ on two sets of equal measure, would
have zero density. In particular, $m$ is even. Also $F(W)>0$ whenever
$W\ne0$. Indeed, $F$ has no nonzero zeros, is positive on nonzero
constant kernels, and continuity along the segment from a nonconstant
$W$ to the constant kernel $1$ rules out a negative value.
Hatami's inequality~\cite[proof of Theorem~2.16]{hatami2010}, which
is proved for possibly asymmetric kernels, gives
\[
    F(U+V)+F(U-V)
    \le 2^{m-1}\bigl(F(U)+F(V)\bigr).
\]
Consequently,
\[
    F\left(\frac{U+V}{2}\right)
    \le\frac{F(U)+F(V)}2-2^{-m}F(U-V),
\]
which implies $F(\frac{U+V}{2})<\frac{F(U)+F(V)}2$ when $U\ne V$. Since $F=\|\cdot\|_H^m$ is convex,
strict midpoint convexity implies strict convexity.

Conversely, suppose that $F$ is strictly convex on real-valued kernels.
Applying strict convexity to $1$ and $-1$ shows that $m$ is even.
Applying it to $W$ and $-W$ then gives $F(W)>0$ for every nonzero
$W$. The same normalization argument
as above, now without taking pointwise absolute values inside $F$,
proves the triangle inequality for $F^{1/m}$. Hence $H$ is norming.
\end{proof}

\section{Forcing Implies OS/OAS}\label{sec:forcing implies OS/OAS}

We first prove Proposition~\ref{prop:forcingandOS}, i.e. symmetry-forcing implies OS or OAS. 

\begin{proof}[Proof of Proposition~\ref{prop:forcingandOS}]
For a fixed directed graph $D$, write $\Delta(W)
    :=
    t(D,W)-t\left(D,\frac{W+W^T}{2}\right)$. 
Suppose for contradiction that $D$ is neither OS nor OAS. Then there
exist digraphons $W_+$ and $W_-$ such that
$\Delta(W_+)>0$ and $\Delta(W_-)<0$.

For each $p\in(0,1)$, partition $[0,1]$ into intervals
$I_p=[0,p]$ and $J_p=(p,1]$. Define a digraphon $W_p$ by placing a
rescaled copy of $W_+$ on $I_p\times I_p$, a rescaled copy of $W_-$ on
$J_p\times J_p$, and setting $W_p=0$ on the two off-diagonal blocks.
More precisely,
\[
W_p(x,y)=
\begin{cases}
W_+\left(\frac{x}{p},\frac{y}{p}\right),
    & x,y\in I_p,\\[2mm]
W_-\left(\frac{x-p}{1-p},\frac{y-p}{1-p}\right),
    & x,y\in J_p,\\[2mm]
0, & \text{otherwise}.
\end{cases}
\]
Thus $W_0$ is a copy of $W_-$ and $W_1$ is a copy of $W_+$, so $\Delta(W_0)<0$ and $\Delta(W_1)>0$.

The function $p\mapsto\Delta(W_p)$ is continuous. Indeed, both $t(D,W_p)$ and
$t\left(D,\frac{W_p+W_p^T}{2}\right)$ are polynomials in $p$. Hence, by the intermediate value
theorem, there exists $p_0\in(0,1)$ such that $\Delta(W_{p_0})=0$, i.e. $t(D,W_{p_0})
    =
    t\left(D,\frac{W_{p_0}+W_{p_0}^T}{2}\right)$.

Note that $W_{p_0}$ is not symmetric since neither $W_+$ nor
$W_-$ is symmetric. This contradicts the assumption that $D$ is
symmetry-forcing. Therefore, $D$ is either OS or OAS.
\end{proof}

\section{OS Directed Graphs}\label{section:OS}

In this section, we characterize connected OS directed graphs, prove the strict-convexity characterization in Theorem~\ref{thm:strict-convexity}, and deduce Theorems~\ref{thm:OS_forcing} and~\ref{thm:cyclic-step-forcing}.

\subsection{Characterization of OS directed graphs}
We first prove that the conditions in Theorem~\ref{thm:charcOS} are sufficient.
\begin{thm}\label{thm:os-sufficient-condition}
Let $D$ be a  directed graph satisfying the following properties:
\begin{itemize}
    \item[(a)] $D$ is homomorphic to a directed edge;
    \item[(b)] $D$ is isomorphic to the digraph $D^T$ obtained by reversing the orientation of each edge of $D$;
    \item[(c)] its underlying graph $\overline{D}$ is weakly norming with respect to a bipartition under which all edges of $D$ go from the first part to the second.
\end{itemize}
Then $D$ is Oriented Sidorenko. 
\end{thm}

\begin{proof}
Let $A\cup B$ be a bipartition of $\overline{D}$ such that all edges in $D$ are oriented from $A$ to $B$. Under this bipartition, for any digraphon $W$ we have  $t(D,W)=t_{bi}(\overline{D},W)$. Let $W$ be a digraphon, and let $S=\frac{W+W^T}{2}$. Since $D\cong D^T$, we have
$$
t_{bi}(\overline{D},W)=t(D,W)=t(D^T,W)=t(D,W^T)=t_{bi}(\overline{D},W^T),
$$
and hence $\|W\|_{r(\overline{D})}=\|W^T\|_{r(\overline{D})}$.
By the triangle inequality of $\|\cdot\|_{r(\overline{D})}$, 
$$
\|S\|_{r(\overline{D})}=\|\frac{W+W^T}{2}\|_{r(\overline{D})}\le\frac{1}{2}\|W\|_{r(\overline{D})}+\frac{1}{2}\|W^T\|_{r(\overline{D})}=\|W\|_{r(\overline{D})}, 
$$
which implies $t(D,S)=t_{bi}(\overline{D},S)\le t_{bi}(\overline{D},W)=t(D,W)$. Therefore, $D$ is OS.
\end{proof}

Next, we prove separately that the conditions in Theorem~\ref{thm:charcOS} are also necessary.

\begin{lem}\label{lem:os imply hom to an edge}
If $D$ is an OS directed graph, then $D$ is homomorphic to a directed edge.
\end{lem}

\begin{proof}
Assume for contradiction that $D$ is not homomorphic to a directed edge. It suffices to find a digraph or digraphon $W$ such that
\begin{equation}\label{eq:contra OS}
t(D, W) < t\left(D,\frac{W+{W}^T}{2}\right).
\end{equation}

Let $G$ be the underlying graph of $D$. If $G$ is bipartite, then letting $W$ be a directed edge would suffice. Thus we may assume that $G$ is not bipartite.

Recall that a directed graph is \emph{acyclic} if it contains no directed cycle, or equivalently, if its vertices can be ordered so that every edge is oriented from left to right. If $D$ is not acyclic, let $W$ be any acyclic orientation of $G$, obtained by fixing an ordering of $V(G)$ and orienting every edge from left to right. Since a directed graph admitting a homomorphism to an acyclic directed graph must itself be acyclic, we have $t(D,W)=0$. However, the identity map of $G$ gives $t\left(D,\frac{W+W^T}{2}\right)>0$. Hence \eqref{eq:contra OS} holds. We may therefore further assume that $D$ is acyclic.

Let $C$ be a minimal subgraph of $G$ such that $G$ is homomorphic to $C$, that is, any proper subgraph of $C$ is not a homomorphic image of $G$. Note that $C$ is not bipartite since $G$ is not. In particular, this implies that $C$ contains a cycle, and thus there exists an orientation $W$ of $C$ that is not acyclic. We claim that this $W$ would satisfy (\ref{eq:contra OS}). 

Indeed, since $G$ is homomorphic to $C$, the right-hand side of (\ref{eq:contra OS}) is positive. On the other hand, let $C'$ be a sub-digraph of $D$ whose underlying graph is $C$. If there exists a homomorphism of $D$ to $W$, then restricting to $C'$ must give an isomorphism from $C'$ to $W$ due to the minimality of the underlying graph $C$. However, note that $C'$ is acyclic and $W$ is not, contradicting the isomorphism. This completes the proof.
\end{proof}

Given a graph $H$, a graphon $W$, two vertices $u,w\in V(H)$ and two numbers $x,y\in [0,1]$, define the rooted density $t(H, W~|~u\rightarrow x, w\rightarrow y)$ to be the density of $H$ in $W$ where $u$ is mapped to $x$ and $w$ is mapped to $y$, i.e.
$$
t(H, W~|~u\rightarrow x, w\rightarrow y)=\int\prod_{e\in E(H)}W(x_{e})\prod_{v\in V(H)\setminus \{u,w\}}dx_v,~\text{where $x_{u}=x$ and $x_{w}=y$}.
$$

\begin{lem}\label{lem:reverse bijection}
Let $D$ be an OS/OAS directed graph with no isolated vertices. There exists a bijection $f:E(D)\rightarrow E(D)$ satisfying the following property: for any $(u,w),(u',w')\in E(D)$ with $(u',w')=f((u,w))$ there exists an automorphism $g$ of $\overline{D}$ such that $g(u)=w'$ and $g(w)=u'$. 
\end{lem}

\begin{proof}
Let $R$ be a double star whose centers have degrees $2v(D)$ and $3v(D)$. More precisely, $R$ is a tree with exactly two non-leaf vertices $r_1$ and $r_2$ such that $d(r_1)=2v(D)$ and $d(r_2)=3v(D)$. For any fixed graphon $S$, we define $U$ such that 
$$
\begin{aligned}
U(x,y)&=\frac{1}{2}(t(R, S|r_1\rightarrow x, r_2\rightarrow y)-t(R, S|r_1\rightarrow y, r_2\rightarrow x))\\
&=\frac{1}{2}S(x,y)(t(R\setminus\{r_1r_2\}, S|r_1\rightarrow x, r_2\rightarrow y)-t(R\setminus\{r_1r_2\}, S|r_1\rightarrow y, r_2\rightarrow x)).
\end{aligned}
$$ 

Note that $U$ is skew-symmetric and $|U(x,y)|\le S(x,y)$. Let $W_\eps=\frac{1}{2}(S+\eps U)$ where $\eps\in[-1,1]$, then $W_{\eps}$ is a digraphon and $\frac{W_{\eps}+W_{\eps}^T}{2}=\frac{S}{2}$. Let $f_S$ be a function on $[-1,1]$ defined by $f_S(\eps)=t(D, W_\eps)$. Note that $f_S$ is a polynomial of $\eps$. Since $D$ is OS/OAS, by definition $f_S(\eps)$ has a local minimum/maximum at $0$, and hence $f_S'(0)=0$, i.e. 
\begin{equation}\label{eq:derivative equal zero}
\sum_{e=(u,w)\in E({D})}\int U(x,y)t(\overline{D\setminus e}, S~|~u\rightarrow x,w\rightarrow y)dxdy=0.
\end{equation}

For any $e=(u,w)\in E(D)$, let $D_e$ be the graph obtained by gluing $\overline{D}$ and $R$ along the edge $(u,w)$ and $(r_1,r_2)$. More precisely,  $V(D_e)=V(D)\cup V(R)\setminus\{r_1,r_2\}$ and 
$$
E(D_e)=E(\overline{D})\cup\{uv~|~r_1v\in E(R),v\in V(R)\setminus\{r_2\}\}\cup \{wv~|~r_2v\in E(R),v\in V(R)\setminus\{r_1\}\}.
$$
Let $D_{e^T}$ be the graph obtained by gluing $\overline{D}$ and $R$ along the $(w,u)$ and $(r_1,r_2)$. That is,  $V(D_{e^T})=V(D)\cup V(R)\setminus\{r_1,r_2\}$ and 
$$
E(D_{e^T})=E(\overline{D})\cup\{wv~|~r_1v\in E(R),v\in V(R)\setminus\{r_2\}\}\cup \{uv~|~r_2v\in E(R),v\in V(R)\setminus\{r_1\}\}.
$$
By definition, (\ref{eq:derivative equal zero}) is equivalent to the following equation.
\begin{equation*}
\sum_{e}t(D_e,S)=\sum_{e}t(D_{e^T},S).
\end{equation*}

Since $D$ has no isolated vertices, neither do the graphs
$D_e$ and $D_{e^T}$. By Lovász's linear independence
result~\cite[Proposition~5.44(c)]{lovasz2012large}, the preceding
identity implies that each isomorphism class occurs with the
same multiplicity in the two multisets
$\{D_e\}_{e\in E(D)}$ and $\{D_{e^T}\}_{e\in E(D)}$.
Consequently, there exists a bijection
$f:E(D)\to E(D)$ such that $D_e\cong D_{f(e)^T}$ for every $e\in E(D)$.

Fix $e=(u,w)$ and write $f(e)=e'=(u',w')$. Let $g:D_e\to D_{e'^T}$ be an isomorphism. Note that $w$ is the only vertex in $D_{e}$ with degree at least $3v(D)$ and $u'$ is the only vertex in $D_{e'^T}$ with degree at least $3v(D)$, so we must have $g(w)=u'$. Similarly, we have $g(u)=w'$. By permuting pendant leaves with the same neighbour if necessary, we may
assume without loss of generality that $g$ maps newly added leaves to
newly added leaves. Restricting $g$ to $V(\overline D)$ then gives an
automorphism of $\overline D$ with $g(u)=w'$ and $g(w)=u'$. Therefore,
$f$ is a bijection with the desired properties.
\end{proof}

\begin{lem}\label{lem:os imply symmetry}
If $D$ is a connected OS directed graph, then $D\cong D^T$.
\end{lem}

\begin{proof}
The edgeless case is immediate, so assume $e(D)>0$.
By Lemma~\ref{lem:os imply hom to an edge}, $D$ is homomorphic to a directed edge. Let $X$ and $Y$ be a bipartition of $V(D)$ such that all edges of $D$ are oriented from $X$ to $Y$. By Lemma~\ref{lem:reverse bijection}, there exists a bijection $f:E(D)\rightarrow E(D)$ such that for any $e=(u,w),e'=(u',w')\in E(D)$ with $f(e)=e'$, there exists an automorphism $g$ of $\overline{D}$ such that $g(u)=w'$ and $g(w)=u'$. Since $\overline{D}$ is connected, any such automorphism maps $X$ to $Y$ and $Y$ to $X$, which implies $D\cong D^T$.
\end{proof}

\begin{lem}\label{lem:os imply weakly norming}
If $D$ is a connected OS directed graph with at least one edge, then $\overline{D}$ is weakly norming.
\end{lem}

\begin{proof}
By Lemmas~\ref{lem:os imply hom to an edge}
and~\ref{lem:os imply symmetry}, there is a bipartition $A\cup B$
of $V(D)$ such that every edge is directed from $A$ to $B$, and
$D\cong D^T$. Put $H=\overline D$ with this fixed bipartition and
$F(U)=t_{bi}(H,U)=t(D,U)$. Then $F(U^T)=F(U)$ for every kernel $U$.

Given two digraphons $U,V$, split the host space into two parts of
measure $1/2$ and define the block kernel
\[
    W=\begin{pmatrix}0&U\\ V^T&0\end{pmatrix}.
\]
Since $D$ is connected, a homomorphism contributing to $t(D,W)$
maps $A$ and $B$ to opposite host parts. Thus, writing $n=v(D)$,
\[
    t(D,W)=2^{-n}\bigl(F(U)+F(V^T)\bigr)
           =2^{-n}\bigl(F(U)+F(V)\bigr).
\]
The off-diagonal blocks of $(W+W^T)/2$ are $(U+V)/2$ and its
transpose, so similarly
\[
    t\left(D,\frac{W+W^T}{2}\right)
    =2^{1-n}F\left(\frac{U+V}{2}\right).
\]
The OS inequality therefore gives
\[
    F\left(\frac{U+V}{2}\right)
    \le\frac{F(U)+F(V)}2.
\]
Scaling extends this inequality to all nonnegative bounded kernels.
For each pair $U,V$, the function $s\mapsto F((1-s)U+sV)$ is a
polynomial and hence continuous. Midpoint convexity therefore implies
convexity. Lemma~\ref{lem:kernel-convexity}(1) now shows directly
that $H$ is weakly norming on the full kernel space.
\end{proof}

By Theorem~\ref{thm:os-sufficient-condition}, Lemma~\ref{lem:os imply hom to an edge}, Lemma~\ref{lem:os imply symmetry}, and Lemma~\ref{lem:os imply weakly norming}, Theorem~\ref{thm:charcOS} follows immediately.

\subsection{Variants of the Sidorenko-forcing implication}\label{sec:Sidorenko-forcing}

We first prove the strict-convexity characterization for arbitrary weakly norming bipartite graphs containing a cycle. We then deduce symmetry forcing for connected OS directed graphs and the corresponding step-forcing result.

\subsubsection{Strict convexity}

The following equality lemma is the main step. Its proof uses the weak
H\"older inequality and an auxiliary sign coordinate.

\begin{lem}\label{lem:sign-lift-equality}
Let $H=(X\cup Y,E)$ be weakly norming and contain a cycle. Suppose
that $S,A$ are kernels such that $S\ge|A|$ and
\[
    t_{bi}(H,S+A)=t_{bi}(H,S-A)=t_{bi}(H,S).
\]
Then $A=0$ almost everywhere.
\end{lem}

\begin{proof}
Write $\Phi(K)=t_{bi}(H,K)$ and $c=\Phi(S)$. On the product
probability space $\Omega=[0,1]\times\{-1,1\}$, define, for each
kernel $B$,
$$
    K_B\bigl((x,\varepsilon),(y,\delta)\bigr)
    =S(x,y)+\varepsilon\delta B(x,y).
$$
Alternatively, $K_B$ can be viewed as a block kernel
$$
\begin{bmatrix}
    S+B& S-B\\
    S-B& S+B
\end{bmatrix}.
$$
Note that $\Phi(K_B)$ can be expressed as the expectation
$$
\mathbb E\left[\prod_{(u,v)\in E(H)}\left(S(x_u,y_v)+\eps_u\delta_vB(x_u,y_v)\right)\right]
$$
where for every $u\in X$ and $v\in Y$ the variables $(x_u,\eps_u)$ and $(y_v,\delta_v)$ are sampled uniformly and independently from $\Omega$. Sampling $(x_u,-\eps_u)$ instead of $(x_u,\eps_u)$ for each $u\in X$ does not change the distribution. Thus, $\Phi(K_B)=\Phi(K_{-B})$. Whenever $|B|\le S$, both kernels
are nonnegative and convexity gives
\begin{equation}\label{eq:sign-lift-lower}
    \Phi(K_B)=\frac{\Phi(K_B)+\Phi(K_{-B})}{2}\ge\Phi\left(\frac{K_B+K_{-B}}2\right)=\Phi(K_0)=c.
\end{equation}
For $B=A$, fix all vertex signs before integrating the base coordinates.
Each edge is then decorated by either $S+A$ or $S-A$. By the weak
H\"older inequality (Lemma~\ref{lem:weak-holder}) and the equality $\Phi(S+A)=\Phi(S-A)=c$, each such decorated density is at most $c$.
Averaging the signs and using~\eqref{eq:sign-lift-lower} gives
$\Phi(K_A)=c$.

Let $r,s:[0,1]\to\{-1,1\}$ be measurable. The separate
measure-preserving changes
$(x,\varepsilon)\mapsto(x,r(x)\varepsilon)$ and
$(y,\delta)\mapsto(y,s(y)\delta)$ give
\begin{equation}\label{eq:sign-switch}
    \Phi(K_{Ars})=c,
    \qquad (Ars)(x,y)=A(x,y)r(x)s(y).
\end{equation}

We claim that
\begin{equation}\label{eq:sign-multiplier}
    \Phi(K_{Ah})=c
\end{equation}
for every bounded measurable function $h:[0,1]^2\to\mathbb R$.

First suppose that
\[
    h(x,y)=\sum_{i=1}^k a_i r_i(x)s_i(y),
    \qquad L=\sum_{i=1}^k|a_i|,
\]
is a finite linear combination of products of sign functions.
If $L=0$, then $h=0$ and $\Phi(K_{Ah})=\Phi(K_0)=c$.
We may therefore assume that $L>0$.
For $0\le t\le1/L$, we have
\[
    K_{tAh}=(1-tL)K_0+
    \sum_{i=1}^k t|a_i|K_{\operatorname{sgn}(a_i)Ar_is_i}.
\]
The right-hand side is a convex combination of nonnegative kernels,
each of density $c$ by~\eqref{eq:sign-switch}. Hence convexity gives
$\Phi(K_{tAh})\le c$. Moreover, $|tAh|\le|A|\le S$, so
\eqref{eq:sign-lift-lower} gives the reverse inequality.
Thus $\Phi(K_{tAh})=c$ for all $0\le t\le1/L$.
Since $t\mapsto\Phi(K_{tAh})$ is a polynomial, this identity holds
for every real $t$. Taking $t=1$ proves~\eqref{eq:sign-multiplier}
for such $h$.

Now let $h$ be any bounded measurable function.
Choose rectangular step functions $h_n$ such that
$\|h_n\|_\infty\le\|h\|_\infty$ and $h_n\to h$ in $L^1$.
Every rectangular step function is a finite linear combination
of products of sign functions, since
$\mathbf1_P=(1+r_P)/2$ for $r_P=2\mathbf1_P-1$.
Thus $\Phi(K_{Ah_n})=c$ for every $n$. Furthermore,
\[
    \|K_{Ah_n}-K_{Ah}\|_{L^1(\Omega^2)}
    =\|A(h_n-h)\|_{L^1([0,1]^2)}
    \le\|A\|_\infty\|h_n-h\|_{L^1([0,1]^2)}
    \longrightarrow0.
\]
These kernels are uniformly bounded in absolute value by
$\|S\|_\infty+\|A\|_\infty\|h\|_\infty$.
Therefore, \eqref{eq:kernel-L1-continuity} gives
\[
    \Phi(K_{Ah})=\lim_{n\to\infty}\Phi(K_{Ah_n})=c,
\]
proving~\eqref{eq:sign-multiplier} for every bounded measurable $h$.
In particular, taking $h=\operatorname{sgn}A$, with
$\operatorname{sgn}0=0$, yields
\begin{equation}\label{eq:sign-absolute}
    \Phi(K_{|A|})=c.
\end{equation}

Call $F\subseteq E$ \emph{even} if every degree in $(V(H),F)$ is
even. Expanding the edge product in $\Phi(K_{|A|})$ and averaging
the independent vertex signs leaves precisely the terms indexed by
even edge sets. Thus
$$
    \Phi(K_{|A|})=
    \sum_{F\ \mathrm{even}}
    t_{bi}\bigl(H;|A|\text{ on }F,\ S\text{ on }E\setminus F\bigr).
$$
Every term is nonnegative. The empty set contributes $\Phi(S)=c$. If $C$ is any cycle of $H$, its edge set also occurs in this sum, and its contribution is at least $\Phi(|A|)$ because $S\ge|A|$. Consequently,
\[
    c=\Phi(K_{|A|})\ge c+\Phi(|A|).
\]
It follows that $\Phi(|A|)=0$. Positive definiteness of
$\|\cdot\|_{r(H)}$ gives $A=0$ almost everywhere.
\end{proof}

\begin{proof}[Proof of Theorem~\ref{thm:strict-convexity}]
Suppose first that $H$ is weakly norming and contains a cycle. Put
$\Phi(W)=t_{bi}(H,W)$ and $m=e(H)>1$. By
Lemma~\ref{lem:kernel-convexity}(1), $\Phi$ is convex on nonnegative
kernels. Let $U,V\ge0$ and set $S=(U+V)/2$, $A=(U-V)/2$.
Suppose that midpoint equality holds. The triangle inequality for
the weak norm gives
\begin{equation}\label{eq:convexity}
    \Phi(S)
    \le\left(\frac{\Phi(U)^{1/m}+\Phi(V)^{1/m}}2\right)^m
    \le\frac{\Phi(U)+\Phi(V)}2.
\end{equation}
Since $z\mapsto z^m$ is strictly convex on $[0,\infty)$, equality
throughout~\eqref{eq:convexity} implies
$\Phi(U)=\Phi(V)=\Phi(S)$. As $S\ge|A|$, Lemma~\ref{lem:sign-lift-equality} gives
$A=0$, and hence $U=V$ almost everywhere. This proves strict midpoint
convexity. For completeness, if $U\ne V$ and $0<\lambda\le1/2$,
convexity gives
\[
\begin{aligned}
    \Phi(\lambda U+(1-\lambda)V)
    &\le2\lambda\Phi((U+V)/2)+(1-2\lambda)\Phi(V)\\
    &<\lambda\Phi(U)+(1-\lambda)\Phi(V).
\end{aligned}
\]
Interchanging $U,V$ treats $1/2\le\lambda<1$.

Conversely, if $\Phi$ is strictly convex on nonnegative kernels, then
it is convex, so Lemma~\ref{lem:kernel-convexity}(1) shows that $H$
is weakly norming. To exclude forests, choose
$f:[0,1]\to\{-1,1\}$ with $\int f=0$, and put
\begin{equation}\label{eq:forest-kernels}
    U(x,y)=\frac12+\frac14f(x)f(y),
    \qquad
    V(x,y)=\frac12-\frac14f(x)f(y).
\end{equation}
These are distinct graphons with all row and column integrals equal
to $1/2$, and $(U+V)/2\equiv1/2$. For any bipartite forest $F$,
successively integrating out leaves gives
\begin{equation}\label{eq:forest-regular-density}
    t_{bi}(F,U)=t_{bi}(F,V)
    =t_{bi}\left(F,\frac{U+V}{2}\right)=2^{-e(F)}.
\end{equation}
Thus a forest cannot have a strictly convex density functional, and $H$ contains a cycle.
\end{proof}

\subsubsection{Symmetry forcing}

\begin{proof}[Proof of Theorem~\ref{thm:OS_forcing}]
Suppose first that $D$ is connected and OS and that $\overline D$
contains a cycle. By Theorem~\ref{thm:charcOS}, $D$ is a one-direction
orientation of a weakly norming bipartite graph $H=\overline D$, and
$D\cong D^T$. With the bipartition induced by the directions,
\[
    t(D,W)=t_{bi}(H,W)=t_{bi}(H,W^T)
\]
for every digraphon $W$. Theorem~\ref{thm:strict-convexity} gives
\[
    t\left(D,\frac{W+W^T}{2}\right)
    \le\frac{t_{bi}(H,W)+t_{bi}(H,W^T)}2=t(D,W),
\]
with equality if and only if $W=W^T$ almost everywhere. Thus $D$
is symmetry-forcing.

Conversely, no orientation of a forest is symmetry-forcing. To see
this, partition $[0,1]$ into three equal parts and take the step
digraphon
\[
    W=\begin{pmatrix}
        1/2&1&0\\
        0&1/2&1\\
        1&0&1/2
    \end{pmatrix}.
\]
Its row and column integrals are all $1/2$, its symmetrization is
the constant kernel $1/2$, and $W\ne W^T$. Integrating out leaves,
regardless of the direction of each incident edge, gives
\[
    t(D,W)=2^{-e(D)}
    =t\left(D,\frac{W+W^T}{2}\right).
\]
The same equality holds when $D$ has no edges. This proves the
necessary direction and completes the characterization.
\end{proof}

\subsubsection{Cut distance identification and step forcing}\label{subsec:Sidorenko-forcing-step}

Here we prove Theorem~\ref{thm:cyclic-step-forcing}.

\begin{proof}[Proof of Theorem~\ref{thm:cyclic-step-forcing}]
As $H$ is weakly norming and contains a cycle, by Theorem~\ref{thm:strict-convexity}, $t_{bi}(H,\cdot)$ is strictly convex on nonnegative kernels. In particular, $t(H,\cdot)$ is strictly convex on symmetric nonnegative kernels.

For any graphon $W$, put $\theta(W):=t(H,W)$. We first claim that for every $\varepsilon>0$ there exists $\eta=\eta_H(\varepsilon)>0$ such that all graphons $U,V$ with $\|U-V\|_\square\geq\varepsilon$ satisfy
\begin{equation}\label{eq:cut-convexity-gap}
    \theta\left(\frac{U+V}{2}\right)
    \leq\frac{\theta(U)+\theta(V)}2-\eta.
\end{equation}
Suppose otherwise. There are sequences of pairs $((U_n,V_n))_{n\ge0}$ with
$\|U_n-V_n\|_\square\geq\varepsilon$ such that
\[
    \frac{\theta(U_n)+\theta(V_n)}2
    -\theta\left(\frac{U_n+V_n}{2}\right)\xrightarrow{n\rightarrow\infty}0.
\]
The usual proof of graphon compactness by weak regularity applies simultaneously to finitely many graphons; see \cite[Chapter~9]{lovasz2012large} for the underlying compactness argument. Consequently, there exist a subsequence $((U_{n_k},V_{n_k}))_{k\ge0}$, measure-preserving bijections $\psi_k$ and graphons $U,V$ such that
\[
    \|U_{n_k}^{\psi_k}-U\|_\square
    +\|V_{n_k}^{\psi_k}-V\|_\square\xrightarrow{k\rightarrow\infty}0.
\]
Since the same bijection is applied to both graphons in each
pair, we retain $\|U-V\|_\square\geq\varepsilon$. Since homomorphism densities are continuous in the cut norm and invariant under measure-preserving bijections, we obtain
\[
    \theta\left(\frac{U+V}{2}\right)
    =\frac{\theta(U)+\theta(V)}2.
\]
Hence, $U=V$ by strict convexity, which contradicts $\|U-V\|_\square\geq\varepsilon$. Thus \eqref{eq:cut-convexity-gap} holds.

We next observe that $\theta$ is weak$^*$ lower semicontinuous on graphons, i.e. if a sequence $(Z_n)_{n\ge 0}$ converges weak$^*$ to $Z$, then $\liminf_{n\to\infty}\theta(Z_n)\ge \theta(Z)$. Indeed, $H$ is weakly norming in the symmetric convention, so it has the step Sidorenko property.
Suppose $Z_n\xrightarrow{n\rightarrow\infty,~w^*}Z$. For every finite measurable partition $\mathcal P$, weak$^*$ convergence gives
$Z_n^{\mathcal P}\to Z^{\mathcal P}$ in $L^1$.
Thus
\[
    \theta(Z^{\mathcal P})
    =\lim_{n\to\infty}\theta(Z_n^{\mathcal P})
    \leq\liminf_{n\to\infty}\theta(Z_n).
\]
Taking $\mathcal P$ to be successively finer dyadic partitions, so that $Z^{\mathcal P}\to Z$ in $L^1$, and using the $L^1$-continuity of homomorphism densities, we obtain
\begin{equation}\label{eq:weak-star-lsc}
    \theta(Z)\leq\liminf_{n\to\infty}\theta(Z_n).
\end{equation}

Now let $U\prec W$. By the definition of the envelope, there are measure-preserving bijections $\phi_n$ such that $W^{\phi_n}\xrightarrow{n\rightarrow\infty,~w^*}U$. Moreover, $\varepsilon:=\delta_\square(U,W)>0$, since equality of envelopes is equivalent to zero cut distance. For every $n$ we have $\|U-W^{\phi_n}\|_\square\geq\varepsilon$. Hence, by \eqref{eq:cut-convexity-gap}, there exists $\eta>0$ such that 
$$
    \theta\left(\frac{U+W^{\phi_n}}2\right)
    \leq\frac{\theta(U)+\theta(W)}2-\eta.
$$
Since $(U+W^{\phi_n})/2\xrightarrow{w^*}U$, by \eqref{eq:weak-star-lsc},
$$
    \theta(U)
    \leq\liminf_{n\to\infty}
        \theta\left(\frac{U+W^{\phi_n}}2\right)
    \leq\frac{\theta(U)+\theta(W)}2-\eta.
$$
Hence $\theta(W)\geq\theta(U)+2\eta>\theta(U)$,
proving that $t(H,\cdot)$ is cut distance identifying.

Finally, cut distance identification implies step forcing: whenever $W\ne W^{\mathcal P}$ as kernels, we have $W^{\mathcal P}\prec W$, and consequently $t(H,W^{\mathcal P})<t(H,W)$; see \cite[Section~3.1]{DolezalGrebikHladkyRochaRozhon2022}. 
\end{proof}

\begin{proof}[Proof of Corollary~\ref{cor:char-step-forcing}]
Let $H$ be connected. If $H$ is step forcing, then it is step Sidorenko
by definition. It cannot be a forest: for the graphon $U$
in~\eqref{eq:forest-kernels} and the one-part partition $\mathcal P$,
we have $U^{\mathcal P}\equiv1/2\ne U$ but
$t(H,U^{\mathcal P})=t(H,U)$ by~\eqref{eq:forest-regular-density}.

Conversely, suppose that $H$ is step Sidorenko and contains a cycle. Dole\v{z}al et al.~\cite{DolezalGrebikHladkyRochaRozhon2022} show that a connected step Sidorenko graph is weakly norming. Theorem~\ref{thm:cyclic-step-forcing} therefore implies that $H$ is step forcing.
\end{proof}

\section{OAS Symmetry-Forcing Directed Graphs}\label{sec:OAS} 

We construct the following family. 

\begin{definition}\label{def:Dt}
For an integer $t\geq 2$, let $D_t$ be
the directed graph whose underlying graph is $K_{t,2^t}$. Let $X=\{u_1,\ldots,u_t\}$ be the part of size $t$, and index the vertices in the other part by the
vectors $\eps=(\eps_1,\ldots,\eps_t)\in\{-1,1\}^t$. Denote the vertex
corresponding to $\eps$ by $v_{\eps}$. For every $i\in[t]$, orient the
edge between $u_i$ and $v_{\eps}$ as
\[
    u_i\longrightarrow v_{\eps}
    \quad\text{if }\eps_i=1,
    \qquad\text{and}\qquad
    v_{\eps}\longrightarrow u_i
    \quad\text{if }\eps_i=-1.
\]    
\end{definition}

Thus, the $2^t$ vertices in the second part realize all possible orientations of their incident edges. Next we prove that all $D_t$ are OAS and symmetry-forcing, i.e. Theorem~\ref{thm:Dt OAS symmetry}.

\begin{proof}[Proof of Theorem~\ref{thm:Dt OAS symmetry}]
Let $W$ be a digraphon, and write $S=\frac{W+W^T}{2}$ and $A=\frac{W-W^T}{2}$. For $\sigma\in\{-1,1\}$, set $W_\sigma=S+\sigma A$, so that
$W_1=W$ and $W_{-1}=W^T$. For
$\eps=(\eps_1,\ldots,\eps_t)\in\{-1,1\}^t$ and
$\boldsymbol{x}=(x_1,\ldots,x_t)\in[0,1]^t$, define $P_\eps(\boldsymbol{x})
=
\int_0^1\prod_{i=1}^t W_{\eps_i}(x_i,y)\,dy$. Note that

\begin{equation}\label{eq:Dt-sf1}
\frac{1}{2^t}\sum_\eps P_\eps(\boldsymbol{x})
=
\int_0^1
\prod_{i=1}^t
\frac{W(x_i,y)+W(y,x_i)}{2}\,dy
=
\int_0^1\prod_{i=1}^t S(x_i,y)\,dy.    
\end{equation}

By definition of $D_t$, AM-GM and \eqref{eq:Dt-sf1},
$$
\begin{aligned}
t(D_t,W)
& =\int_{[0,1]^t}\prod_{\eps\in\{-1,1\}^t}P_\eps(\boldsymbol{x})\,d\boldsymbol{x}
\le \int_{[0,1]^t}\left(
\frac{1}{2^t}
\sum_{\eps\in\{-1,1\}^t}P_\eps(\boldsymbol{x})
\right)^{2^t}\,d\boldsymbol{x}\\
&= \int_{[0,1]^t}\left(
\int_0^1\prod_{i=1}^t S(x_i,y)\,dy
\right)^{2^t}\,d\boldsymbol{x} =t(D_t,S).
\end{aligned}
$$
Hence $D_t$ is OAS.

It remains to prove that $D_t$ is symmetry-forcing. Suppose that $t(D_t,W)=t(D_t,S)$. Then, for almost every $\boldsymbol{x}$, equality holds in the AM--GM inequality above. Thus $P_\eps(\boldsymbol{x})$ is independent of $\eps$ for almost every $\boldsymbol{x}$.

Then 
\begin{equation}\label{eq:Dt-sf2}
\int_0^1
A(x_1,y)A(x_2,y)
\prod_{i=3}^t S(x_i,y)\,dy=\frac{1}{2^t}
\sum_{\eps\in\{-1,1\}^t}
\eps_1\eps_2P_\eps(\boldsymbol{x})=0    
\end{equation}

for almost every $x_1,\ldots,x_t$, where the empty product is $1$ when $t=2$.

Let $f(x)$ be any bounded measurable function. Multiplying \eqref{eq:Dt-sf2} by $f(x_1)f(x_2)$ and then integrating over $x_1,\dots,x_t$ gives
\[
\int_{0}^1\left(\int_0^1 f(x) A(x,y)\,dx\right)^2\left(\int_0^1S(x,y)\,dx\right)^{t-2}\,dy=0.
\]
When $t=2$, the final weight in the integrand is interpreted as $1$.
Take $f$ from the countable family of indicators of dyadic intervals.
There is then a single null set outside which the displayed identity
forces
\[
    \int_0^1 f(x)A(x,y)\,dx=0
\]
for every such $f$ whenever $\int_0^1S(x,y)\,dx>0$.
Dyadic step functions are dense in $L^2([0,1])$, so
$A(\cdot,y)=0$ almost everywhere for each such $y$ outside that null
set. If $\int_0^1S(x,y)\,dx=0$, the inequalities $0\le|A|\le S$
give the same conclusion. Thus $A=0$ almost everywhere and
$W=W^T$ almost everywhere. Therefore, $D_t$ is symmetry-forcing.
\end{proof}

We observe the following interesting example of an OAS graph whose underlying graph is not bipartite. Let $\widetilde{D}_2$ be the directed graph obtained from $D_2$ by adding the directed edge $u_1\to u_2$.

\begin{proof}[Proof of Proposition~\ref{prop:nonbipartite-oas}]
Let $W$ be a digraphon and write $S=\frac{W+W^T}{2}$. As before, let $W_1=W$, $W_{-1}=W^T$ and let $\eps\in \{-1,1\}^2$, $P_{\eps}(x_1,x_2)=\int_0^1W_{\eps_1}(x_1,y)W_{\eps_2}(x_2,y)dy$. Further, let $R(x_1,x_2)=\prod_{\eps\in \{-1,1\}^2}P_{\eps}(x_1,x_2)$. Note that $R$ is symmetric, i.e. $R(x_1,x_2)=R(x_2,x_1)$. Hence,

$$
\begin{aligned}
t(\widetilde{D}_2,W) &= \int_{[0,1]^2}W(x_1,x_2)R(x_1,x_2)\,dx_1\,dx_2 =\frac{1}{2}\int_{[0,1]^2}(W(x_1,x_2)+W(x_2,x_1))R(x_1,x_2)\,dx_1\,dx_2\\
&=\int_{[0,1]^2}S(x_1,x_2)R(x_1,x_2)\,dx_1\,dx_2
\end{aligned}
$$
By AM-GM and \eqref{eq:Dt-sf1}, $R(x_1,x_2)\le \left(\frac{1}{4}\sum_{\eps}P_\eps(x_1,x_2)\right)^4=\left(\int_0^1S(x_1,y)S(x_2,y)dy\right)^4$.
Thus
$$
t(\widetilde{D}_2,W)\le \int_{[0,1]^2}S(x_1,x_2)\left(\int_0^1S(x_1,y)S(x_2,y)dy\right)^4\,dx_1dx_2=t(\widetilde{D}_2,S).
$$
Therefore, $\widetilde{D}_2$ is OAS. Its underlying graph contains the triangle $u_1u_2v_{(1,1)}u_1$, so it is not bipartite. 
\end{proof}

\section{Characterizing stars, paths, and cycles}\label{sec:stars, paths, and cycles}
In this section, we prove Theorems~\ref{thm:OAS stars},~\ref{thm:OAS paths}, and~\ref{thm:OAS cycles}.

\subsection{Stars}
Here we prove Theorem~\ref{thm:OAS stars} which characterize OS and OAS oriented stars.

\begin{proof}[Proof of Theorem~\ref{thm:OAS stars}]
If $D$ has one edge, then for every digraphon $W$,
\[
    t(D,W)=\int_{[0,1]^2}W(x,y)\,dx\,dy
    =\int_{[0,1]^2}\frac{W(x,y)+W(y,x)}2\,dx\,dy
    =t\left(D,\frac{W+W^T}{2}\right).
\]
Thus $D$ is both OS and OAS.

We first prove that if the oriented star $D$ has at least two edges then it is not OS. Let $c$ be its center.
Suppose for contradiction that $D$ is OS. By Theorem~\ref{thm:charcOS}, the orientation $D$ must be a one-direction orientation and must satisfy $ D\cong D^T.$
Since the two parts of the underlying star are the center and the leaves, being a one-direction orientation means that either every edge is directed away from $c$, or every edge is directed toward $c$. The center is the unique vertex of degree greater than one. Therefore,
every isomorphism from $D$ to $D^T$ must map $c$ to itself. However, the center is a source in one of $D$ and $D^T$ and a sink in the other. Hence $D$ cannot be isomorphic to $D^T$, a contradiction. Thus no oriented star with at least two edges is OS.

We now find all oriented stars $D$ which are OAS. Suppose $D$ has at least two edges, with center $c$. First suppose that $c$ has in-degree and out-degree $k$. For a digraphon $W$, define $W^+(x)=\int_0^1W(x,y)dy$ and $W^-(x)=\int_0^1W(y,x)dy$. Conditioning on the image of $c$ gives $t(D,W)=\int_0^1\bigl(W^+(x)W^-(x)\bigr)^kdx$. Let $S=(W+W^T)/2$, then the AM--GM inequality yields
\[
\begin{aligned}
    t(D,W)
    \leq
    \int_0^1
    \left(\frac{W^+(x)+W^-(x)}2\right)^{2k}dx
    =\int_0^1
    \left(\int_0^1S(x,y)\,dy\right)^{2k}dx
    =t(D,S).
\end{aligned}
\]
Thus $D$ is OAS.

Conversely, suppose that $D$ is OAS. Since the center is the unique non-leaf vertex, every automorphism of the underlying star fixes it. By Lemma~\ref{lem:reverse bijection}, the resulting edge bijection maps each edge directed away from the center to one directed toward the center, and vice versa. Hence the two degrees are equal.
\end{proof}

\subsection{Paths}\label{sec:paths}
Here we prove Theorem~\ref{thm:OAS paths} which characterize OS and OAS oriented paths.

\begin{proof}[Proof of Theorem~\ref{thm:OAS paths}]
A directed edge is both OS and OAS by Theorem~\ref{thm:OAS stars}. No orientation of a two-edge path is OS, again by Theorem~\ref{thm:OAS stars}. For a path with at least three edges, its underlying graph is not biregular, whereas every weakly norming bipartite graph is biregular~\cite[Theorem~2.10(ii)]{hatami2010}. Thus Theorem~\ref{thm:charcOS} excludes all such paths from OS.

We now consider the OAS assertion. Every orientation of a two-edge path is an oriented star, and by Theorem~\ref{thm:OAS stars}, an oriented star with at least two edges is OAS if and only if its center has equal in-degree and out-degree. Thus the directed two-edge path is the only OAS orientation of a two-edge path.

It remains to show that no orientation of a path with at least three edges is OAS. Let $P$ be an oriented path with vertices $0,1,\ldots,n$, where $n\geq3$ and $\{i-1,i\}$ is an edge for $1\leq i\leq n$. Suppose for contradiction that $P$ is OAS. The underlying path has only two automorphisms: the identity and the reflection $i\mapsto n-i$. The identity cannot satisfy the endpoint-reversal condition in Lemma~\ref{lem:reverse bijection}. Thus the automorphism must reverse the end points of every edge, which gives 
\begin{equation}\label{eq:path-palindrome}
    (i-1,i)\in E(P)
    \quad\Longleftrightarrow\quad
    (n-i,n-i+1)\in E(P).
\end{equation}
The OAS property is invariant under reversing every edge. Replacing $P$ by $P^T$ if necessary, we may therefore assume that $(0,1)\in E(P)$. There are two cases, according to the orientation of the next edge.

Suppose first that $(1,2)\in E(P)$. By~\eqref{eq:path-palindrome}, the last two edges are $(n-2,n-1)$ and $(n-1,n)$. Partition $[0,1]$ into three sets $X,Y,Z$ of measure $1/3$, and let $W$ be the step digraphon whose matrix with respect to this partition is
$$
    \begin{bmatrix}
        1&\eps^3&0\\
        0&0&\eps^2\\
        0&0&0
    \end{bmatrix}.
$$
The terms of degree less than six in $\eps$ can be read from the possible block assignments of the path. They give
\begin{equation}\label{eq:path-case1-W}
    t(P,W)
    =3^{-(n+1)}\bigl(1+\eps^3+\eps^5\bigr)+O(\eps^6).
\end{equation}
The constant term maps every vertex to $X$. The $\eps^3$ term maps only the final vertex to $Y$, and the $\eps^5$ term maps the last two vertices to $Y,Z$, respectively.

For $S=(W+W^T)/2$, the two ends of the path contribute symmetrically. There are two assignments of weight $\eps^3/2$ and two assignments of weight $\eps^5/4$. Hence
\begin{equation}\label{eq:path-case1-S}
    t(P,S)
    =3^{-(n+1)}\left(1+\eps^3+\frac12\eps^5\right)
      +O(\eps^6).
\end{equation}
It follows that $t(P,W)>t(P,S)$ for all sufficiently small $\eps>0$, a contradiction.

It remains to consider $(2,1)\in E(P)$. In this case \eqref{eq:path-palindrome} gives $(n-1,n-2),(n-1,n)\in E(P)$. Use instead the step digraphon
$$
    W=
    \begin{bmatrix}
        1&0&0\\
        \eps^3&0&\eps^2\\
        0&0&0
    \end{bmatrix}.
$$
The same block-assignment calculation gives
$$
\begin{aligned}
    t(P,W)
    &=3^{-(n+1)}\bigl(1+\eps^3+\eps^5\bigr)+O(\eps^6),\\
    t(P,S)
    &=3^{-(n+1)}\left(1+\eps^3+\frac12\eps^5\right)
      +O(\eps^6).
\end{aligned}
$$
Again $t(P,W)>t(P,S)$ for sufficiently small $\eps>0$. This completes the proof.
\end{proof}

\subsection{Cycles}
The following lemma is the main tool for characterizing OAS cycles.

We first prove a necessary condition for a general OAS graph by expanding its density around a constant digraphon. Such local
expansions also appear in the study of local Sidorenko inequalities; see Lov\'asz~\cite{Lovasz2011local} and
Fox--Wei~\cite{FoxWei2017local}. Here the additional phase in our Fourier perturbation records the orientation imbalance of the
shortest cycles.

Let $D$ be a directed graph, and let $g$ be the girth of its underlying graph $\overline D$. Denote by $\mathcal C_g(D)$ the set of cycles of length $g$ in $\overline D$.

For each $C\in\mathcal C_g(D)$, choose one of the two cyclic directions around $C$. Let $f_D(C)$ be the number of edges of $C$ whose orientation in $D$ agrees with this cyclic direction, and let $b_D(C)$ be the number whose orientation disagrees with it. Define the \emph{orientation imbalance} of $C$ by     $\omega_D(C):=f_D(C)-b_D(C).$ Reversing the chosen cyclic direction replaces $\omega_D(C)$ by $-\omega_D(C)$. Thus $\cos\bigl(\omega_D(C)\theta\bigr)$ does not depend on the choice of cyclic direction.

\begin{lem}
\label{lem:girthcondition}
Let $D$ be a directed graph whose underlying graph contains a cycle, and let $g$ be the girth of $\overline D$. If $D$ is OAS, then
\[
    \sum_{C\in\mathcal C_g(D)}
    \cos\bigl(\omega_D(C)\theta\bigr)
    \leq
    |\mathcal C_g(D)|\cos^g(\theta)
\]
for every $\theta\in[0,2\pi)$.
\end{lem}

\begin{proof}
Fix $\theta\in[0,2\pi)$ and let $0<\varepsilon<1/2$. Define
\[
    W_{\varepsilon,\theta}(x,y)
    =
    \frac12+\varepsilon
    \cos\bigl(2\pi(x-y)+\theta\bigr)= 
    \frac12\left(
        1
        +\varepsilon e^{\mathrm i\theta}
            e^{2\pi\mathrm i(x-y)}
        +\varepsilon e^{-\mathrm i\theta}
            e^{-2\pi\mathrm i(x-y)}
    \right).
\]
The bound $\varepsilon<1/2$ ensures that $W_{\varepsilon,\theta}$ takes values in $[0,1]$.

Write $m=e(D)$. By definition,
\begin{equation}
    t(D, W_{\eps, \theta}) = 2^{-e(D)}\int \prod_{(u,w)\in E(D)} \left(1+\eps e^{(2\pi(x_u-x_w)+\theta)\mathrm i}+\eps e^{-(2\pi(x_u-x_w)+\theta)\mathrm i} \right)\prod_{v\in V(D)}dx_v.
\end{equation} 
Expanding the product in the definition of $t(D,W_{\varepsilon,\theta})$, we obtain
\[
\begin{aligned}
t(D,W_{\varepsilon,\theta})
&=
2^{-m}
\sum_{\sigma:E(D)\to\{-1,0,1\}}
\varepsilon^{|\operatorname{supp}\sigma|}
e^{\mathrm i\theta\sum_{e\in E(D)}\sigma_e} \\
&\hspace{2cm}\times
\int_{[0,1]^{V(D)}}
\exp\left(
    2\pi\mathrm i
    \sum_{e=(u,v)\in E(D)}
    \sigma_e(x_u-x_v)
\right)
\prod_{v\in V(D)}\mathrm dx_v .
\end{aligned}
\]
Here $\sigma_e=0$ means that the constant term $1$ is chosen for the edge $e = (u,w)\in E(D)$, while $\sigma_e=1$ or $-1$ means that one of the two Fourier terms is chosen.

By orthogonality of the functions
$x\mapsto e^{2\pi\mathrm i kx}$, the integral is nonzero precisely when
\[
    \sum_{\substack{e=(v,w)\\e\in E(D)}}\sigma_e
    =
    \sum_{\substack{e=(u,v)\\e\in E(D)}}\sigma_e
\]
at every vertex $v$. In other words, the selected signed edges must form a circulation. When this condition holds, the integral equals $1$.

Every nonempty circulation contains a cycle in its support. Since the girth of $\overline D$ is $g$, there are no nonconstant terms of degree less than $g$ in $\varepsilon$. A circulation with support of size exactly $g$ must be supported on a cycle $C\in\mathcal C_g(D)$. For each such cycle there are exactly two possibilities, corresponding to the two directions around $C$.

For the circulation following the chosen cyclic direction of $C$, an edge contributes $\sigma_e=1$ if its orientation in $D$ agrees with that direction and $\sigma_e=-1$ otherwise. Hence
\[
    \sum_{e\in C}\sigma_e
    =
    f_D(C)-b_D(C)
    =
    \omega_D(C).
\]
The circulation in the opposite direction contributes $-\omega_D(C)$. The two corresponding phases therefore add up to
\[
    e^{\mathrm i\omega_D(C)\theta}
    +
    e^{-\mathrm i\omega_D(C)\theta}
    =
    2\cos\bigl(\omega_D(C)\theta\bigr).
\]
It follows that
\begin{equation}
\label{eq:fourier-expansion-W}
t(D,W_{\varepsilon,\theta})
=
2^{-m}
\left(
    1
    +
    2\varepsilon^g
    \sum_{C\in\mathcal C_g(D)}
    \cos\bigl(\omega_D(C)\theta\bigr)
    +
    O_D(\varepsilon^{g+1})
\right).
\end{equation}

Now we can write the symmetrization  $S_{\varepsilon,\theta}
    =
    \frac{
        W_{\varepsilon,\theta}
        +
        W_{\varepsilon,\theta}^T
    }2$.
Therefore we have 
\[
\begin{aligned}
S_{\varepsilon,\theta}(x,y)
&=
\frac12
+
\frac{\varepsilon}{2}
\left(
    \cos\bigl(2\pi(x-y)+\theta\bigr)
    +
    \cos\bigl(2\pi(y-x)+\theta\bigr)
\right) \\
&=
\frac12
+
\varepsilon\cos\theta\cos\bigl(2\pi(x-y)\bigr) = 
    \frac12\left(
        1
        +
        \varepsilon\cos\theta\,
            e^{2\pi\mathrm i(x-y)}
        +
        \varepsilon\cos\theta\,
            e^{-2\pi\mathrm i(x-y)}
    \right).
\end{aligned}
\]

The same circulation expansion now gives
\begin{equation}
\label{eq:fourier-expansion-S}
t(D,S_{\varepsilon,\theta})
=
2^{-m}
\left(
    1
    +
    2\varepsilon^g
    |\mathcal C_g(D)|\cos^g\theta
    +
    O_D(\varepsilon^{g+1})
\right).
\end{equation}
Indeed, each of the $g$ selected edges contributes a factor $\cos\theta$, and each shortest cycle again has two possible circulations.

Since $D$ is OAS, $t(D,W_{\varepsilon,\theta})
    \leq
    t(D,S_{\varepsilon,\theta})$.
Substituting \eqref{eq:fourier-expansion-W} and
\eqref{eq:fourier-expansion-S}, we obtain
\[
\begin{aligned}
0
&\leq
2^{1-m}\varepsilon^g
\left(
    |\mathcal C_g(D)|\cos^g\theta
    -
    \sum_{C\in\mathcal C_g(D)}
    \cos\bigl(\omega_D(C)\theta\bigr)
\right)
+
O_D(\varepsilon^{g+1}).
\end{aligned}
\]
Divide by $2^{1-m}\varepsilon^g$ and let $\varepsilon\to0$. This gives $\sum_{C\in\mathcal C_g(D)}
    \cos\bigl(\omega_D(C)\theta\bigr)
    \leq
    |\mathcal C_g(D)|\cos^g\theta,$ as required.
\end{proof}

We now prove that $O_1$ is the only OAS cycle.

\begin{proof}[Proof of Theorem~\ref{thm:OAS cycles}]
The OS assertion follows from Theorem~\ref{thm:charcOS}: its
one-direction condition forces an alternating even cycle, and every such orientation is self-converse with a weakly norming underlying graph. For the OAS assertion, Griffiths proved that $O_1$ is OAS, so it remains to prove necessity.

Let $C$ be an OAS orientation of a cycle of length $g\geq3$. Choose one of the two cyclic directions. After reversing this choice if necessary, we may assume that $t\leq g/2$ edges point against the chosen direction and the remaining $g-t$ edges point with it. Write $\omega:=g-2t.$ Thus $\omega \geq0$ is the absolute value of the orientation imbalance of
$C$.

The underlying graph of $C$ has exactly one cycle of length $g$. Lemma~\ref{lem:girthcondition} therefore gives
\begin{equation}
\label{eq:cycle-condition}
    \cos(\omega\theta)\leq\cos^g\theta
\end{equation}
for every $\theta\in[0,2\pi)$.

We first show that $\omega\leq2.$ Suppose instead that $\omega\geq3$, and take
$\theta=\frac{2\pi}{\omega}.$ Then the left-hand side of \eqref{eq:cycle-condition} is equal to $1$. On the other hand, $ \left|\cos\left(\frac{2\pi}{\omega}\right)\right|<1,$ and hence $ \cos^g\left(\frac{2\pi}{\omega}\right)<1.$ This contradicts \eqref{eq:cycle-condition}. Therefore we have shown \begin{equation}
    \omega \leq 2. \label{eq:omega2}
\end{equation} 

We now use $\omega \leq 2$ to upper bound the cycle length $g$. Consider \eqref{eq:cycle-condition} for $\theta$ close to zero. We have $\cos(\omega\theta) =  1-\frac{\omega^2}{2}\theta^2+O(\theta^4)$
and
$   \cos^g\theta
    =
    1-\frac{g}{2}\theta^2+O(\theta^4).$
Therefore, $0
    \geq
    \cos(\omega\theta)-\cos^g\theta
    =
    \frac{g-\omega^2}{2}\theta^2+O(\theta^4)$.
Letting $\theta\to 0$ gives 
\begin{equation}
    \omega^2\geq g. \label{eq:omegag}
\end{equation}
Combining this with $\omega \leq 2$, we have $g\leq4.$ So it suffices to consider cycles of length at most $4$.

If $g=3$, then $\omega=g-2t$ is odd. Since $\omega\leq2$, we must have $\omega=1$, but then $\omega^2<g$, a contradiction with (\ref{eq:omegag}).

Hence $g=4$. Since $\omega^2\geq4$ and $\omega\leq2$ by (\ref{eq:omegag}) and (\ref{eq:omega2}), we have $\omega=2$. It follows that $t=\frac{g-\omega}{2}=1.$
Thus exactly one edge points against the other three around the cycle. Up to isomorphism, this orientation is $O_1$.
\end{proof}

Combining Theorem~\ref{thm:OAS cycles} with Theorem~\ref{thm:charcOS}, we obtain a complete classification of oriented cycles satisfying either inequality. An oriented cycle is OS if and only if it is an alternating orientation of an even cycle, and it is OAS if and only if it is isomorphic to $O_1$.

\section*{AI Disclosure}
The authors first used OpenAI GPT-5.6 Sol to prove weaker versions of Theorem~\ref{thm:OS_forcing} and Theorem~\ref{thm:cyclic-step-forcing} concerning reflection graphs. After reading the authors realize both these weaker versions stem from a weaker relection-graph version of Theorem~\ref{thm:strict-convexity}. Then the author used GPT-6 Astra to prove the current version of Theorem~\ref{thm:strict-convexity}, and then use it to prove the current version of Theorem~\ref{thm:OS_forcing} and Theorem~\ref{thm:cyclic-step-forcing}. Theorem GPT-5.6 Sol found the constructions in Theorem~\ref{thm:Dt OAS symmetry}. AI was also used to proofread and polish the manuscript. All AI-generated arguments were simplified, rewritten, and independently verified by the authors, who take full responsibility for its correctness.

\bibliographystyle{abbrv}
\bibliography{refs1}

\end{document}